\documentclass[12pt]{article}
\usepackage[left=2.5cm,right=2.5cm,top=3cm,bottom=4cm]{geometry}
\usepackage{authblk} 
\usepackage{verbatim}
\usepackage{framed}
\usepackage{booktabs} 
\usepackage{floatrow}
\newfloatcommand{capbtabbox}{table}[][\FBwidth]
\usepackage{amsmath,amssymb,amsfonts,latexsym,cancel}
\usepackage{amsthm}
\usepackage{relsize,exscale}
\usepackage{enumerate}
\usepackage{color}
\usepackage{xcolor}
\usepackage{graphics}
\usepackage{float}
\usepackage{array}
\usepackage{etoolbox}
\usepackage{dsfont}
\usepackage{amsbsy}
\usepackage{tikz}
\usetikzlibrary{babel}
\usetikzlibrary{arrows,chains,matrix,positioning,scopes}

\makeatletter
\tikzset{join/.code=\tikzset{after node path={%
\ifx\tikzchainprevious\pgfutil@empty\else(\tikzchainprevious)%
edge[every join]#1(\tikzchaincurrent)\fi}}}
\makeatother
\tikzset{>=stealth',every on chain/.append style={join},
         every join/.style={->}}
\tikzstyle{labeled}=[execute at begin node=$\scriptstyle,
   execute at end node=$]

\newtheorem{theorem}{Theorem}[section]

\newtheorem{corollary}[theorem]{Corollary}

\newtheorem{example}[theorem]{Example}

\graphicspath{{graficos/}} 
\usepackage{subfiles} 
\usepackage[style=authoryear,backend=bibtex, isbn=false, doi=false, natbib,dashed=false]{biblatex}
\DeclareNameAlias{sortname}{last-first} 
\renewbibmacro*{volume+number+eid}{%
  \printfield{volume}%
  \setunit*{\addnbthinspace}
  \printfield{number}%
  \setunit{\addcomma\space}%
  \printfield{eid}}
\DeclareFieldFormat[article]{number}{\mkbibparens{#1}}

\title{Ruin probabilities as recurrence sequences in a discrete-time risk process}
\author[1]{David J. Santana}
\author[2]{Luis Rinc\' on}
\affil[1]{Divisi\'on Acad\'emica de Ciencias B\'asicas\\UJAT\\M\'exico}
\affil[2]{Departamento de Matem\'aticas\\Facultad de Ciencias\\UNAM\\M\'exico}

\begin{document}
\date{}
\maketitle

\begin{abstract}
We apply the theory of linear recurrence sequences to find an expression for the  ultimate
ruin probability in a discrete-time risk process. We assume the claims follow an arbitrary
distribution with support $\{0,1,\ldots,m\}$, for some integer $m\ge2$.
The method requires to find the zeroes of an $m$ degree polynomial and to solve a  system of $m$ linear
equations. An approximation is derived from the exact ruin formula and several numerical results
and plots are provided as examples.
{\bf Corresponding author} Luis Rinc\'on lars@ciencias.unam.mx\\
\end{abstract}

\section{Introduction}
The Gerber-Dickson risk process~(\cite{2017dickson}) is given by
\begin{equation}
\label{Gerber-Dickson-model}
U (t) = u+t-\sum_{j= 1}^t Y_j,\quad\mbox{for }\ t=0,1,\ldots
\end{equation}
This is a simple stochastic process that represents the evolution in time of the capital of an insurance company.
Here $U(0)=u\ge0$ is an integer that stands for the initial capital, the insurance company receives
one unit of currency as premium in each unit time period and $Y_1, Y_2, \ldots $ are independent,
identically distributed discrete random variables taking values on $\{0,1,\ldots\}$. These are the
random claims amounts payable  at the end of each period and we denote any of them by $Y$.
The probability function of the
claims is $f(y)={\mathds P}(Y=y)$ and the distribution function is $F(y)=P(Y\le y)$, for
$y= 0,1,\ldots$ We will assume that the
mean value of the claims, denoted by $\mu_Y$, satisfies the net profit condition, that is,
$\mu_Y<1$. When this restriction is not set, in the long run 
the capital reaches values below zero with probability $1$, which is not desirable. Observe
that the net profit condition implies that $f(0)>0$.\\

The time of ruin is defined as the stopping time $\tau=\min\,\{t \ge 1: U (t) \le 0 \}$,
where the minimum of the empty set is defined as infinity.
Given an arbitrary distribution for the claims, a central problem in the theory
of ruin is to find the probability of ultimate ruin (infinite time horizon), which is defined by
\begin{equation*}
\psi(u)={\mathds P}(\tau <\infty \mid U (0) = u).\\
\end{equation*}

The risk process~\eqref{Gerber-Dickson-model} is rather elementary and several
generalizations have been proposed. For example,
let $\{N(t): t=0,1,\ldots\}$ be a counting process such that $N(t)\sim\mbox{Binomial}(t,p)$.
The compound binomial risk process is defined as
\begin{equation}
\label{Compound-binomial-process}
U (t) = u+t-\sum_{j= 1}^{N(t)} X_j,\quad\mbox{for }\ t=0,1,\ldots
\end{equation}
where now claims $X_j$ take only positive values $1,2,\ldots$ and they are independent and identically
distributed as before. The risk process~\eqref{Compound-binomial-process} was introduced
and studied by H. U. Gerber in~(\cite{1988gerber}) and is not really
a generalization of~\eqref{Gerber-Dickson-model} since it can be shown
(see~\cite{Santana-Rincon-2020})
that models~\eqref{Gerber-Dickson-model} and~\eqref{Compound-binomial-process} are equivalent
in the sense that $U(t)$ has the same distribution in both models when the claims
are related by $Y_j=R_j\,X_j$, where $R_1, R_2,\ldots$ is a sequence of i.i.d. $\mbox{Bernoulli}(p)$
random variables. As an example of generalization
of~\eqref{Gerber-Dickson-model}, see~\textcite{2009Trufin-Loisel}  where authors consider a
more elaborate discrete-time risk process where experience rating is taken into account.
In particular,
they derive the logarithmic asymptotic behavior of the ultimate ruin probability for large values
of $u$. The probability of ruin in finite horizon where the reserve is invested in a risky asset
is studied in~\textcite{2003Tang-Tsitsiashvili}.  For some other interesting extensions of the
process~\eqref{Gerber-Dickson-model} see also~\textcite{2018Liu-Wang-Guo},
\textcite{2011Peng-Huang-Wang},
\textcite{cai_2002},
\textcite{cossette_marceau_maume-deschamps_2010},
\textcite{diasparra_romera_2009},
\textcite{Jasiulewicz-Kordecki-2015},
\textcite{sun_yang_2003},
\textcite{wu_chen_guo_jin_2015},
\textcite{YANG-etal200921}.
A comprehensive survey of results on several discrete-time risk
processes can be found in~(\cite{2009Li-Lu-Garrido}).\\

The basic question of finding the ruin probability $\psi(u)$ for our
model~\eqref{Gerber-Dickson-model} and its generalizations remains a hard problem
and in this work we give a hint at why this is so.
The technique that we here present to solve the problem makes use of the theory
of linear recurrence sequences. We have applied this technique to both discrete and continuous time
models in~(\cite{Rincon-Santana-2021, Santana-Rincon-4, Rincon-Santana-2022}).
In those works, however, the recurrence sequences arise from a specific distribution assumed
for the claims. The formulas found are thus specialized to that claims distribution.
On the contrary, in the present work we will see that the structure of the
discrete-time risk process allows to directly identify a recurrence sequence
for the ruin probabilities themselves. The claims distribution is only asked to have finite support
but else arbitrary. Thus, we will see that
solving the recurrence sequence will yield a new formula for the probability of ruin.\\

Before presenting our results, the reader should be aware that there exist in the actuarial literature
several formulas for the ruin probability for discrete-time risk processes. Particularly, for the
compound binomial risk process~\eqref{Compound-binomial-process},
we have the H. U. Gerber's formula~(\cite{1988gerber}). When the definition of time of ruin is
modified so that ruin occurs only when the risk process is strictly negative, we have the E. Shiu's
formula~(\cite{1989shiu}). We also have a formula by S. Li and J. Garrido~(\cite{2002Li}).
All those formulas were obtained by different methods and are rather elaborate as they
are expressed as infinite series involving the partial sums $S_k=X_1+\cdots+X_k$.
These expressions can also
be consulted in the excellent survey~(\cite{2009Li-Lu-Garrido}).

\section{Main result}

We here present the procedure that leads to a new expression for the ultimate ruin probability
for the risk process~\eqref{Gerber-Dickson-model}.
Conditioning on the outcome of the first claim (method also known as first step analysis),
it is easy to show (see~\cite{2017dickson})
that the ruin probability $\psi(u)$ for model~(\ref{Gerber-Dickson-model}) satisfies the
recursive relation
\begin{equation}
\label{Recursive-equation-1}
\psi(u)=\sum\limits_{k=0}^u \psi(u+1-k)f(k) + \overline{F}(u),
\quad\mbox{for } u\ge0,
\end{equation}
where $\overline{F}(u)=1-F(u)$. There are several equivalent forms the
relation~\eqref{Recursive-equation-1} can be written, another version is the equation
\begin{equation}
\label{Recursive-equation-2}
\psi(u)=\sum\limits_{k=0}^{u-1} \psi(u-k)\overline{F}(k)+\sum\limits_{k=u}^{\infty}\overline{F}(k),
\quad\mbox{for } u\ge0.
\end{equation}
Evaluating~\eqref{Recursive-equation-2} at $u=0$ (the empty sum is defined as null) yields
the well known result
$\psi(0)=\mu_Y=\sum_{k=0}^\infty \overline{F}(k)$. Thus,
in the ensuing calculations we will assume that the value of $\psi(0)$ is known to be the average value of
the claims.\\

Let $m\ge2$ be an integer and suppose that the claims take
only the first $m+1$ values $0,1,\ldots,m$ with $f(m)$ strictly positive.
Then we have $\overline{F}(u)=0$ for $u\ge m$ and the
relation~\eqref{Recursive-equation-1} takes the simpler form
\begin{equation}
\label{Recursive-relation-1}
\psi(u)=\sum\limits_{k=0}^m \psi (u+1-k)f(k),\quad\mbox{for } u\ge m.
\end{equation}
Observe that the upper limit of the sum is now the parameter $m$ not $u-1$.
Solving for $\psi (u+1)$ gives
\begin{equation}
\label{Recursive-relation-2}
\psi (u+1)=\frac{1}{f(0)}\,[\ \psi(u)[1-f(1)]-\psi (u-1)f(2)-\cdots-\psi (u-(m-1))f(m)\ ],
\end{equation}
for  $u\ge m$. This is a linear recurrence sequence of order $m$ for $\psi(u)$ with initial
data $\psi(1),\ldots,\psi(m)$. These initial values of the sequence are given by the first terms
of~\eqref{Recursive-equation-1}, namely,
\begin{eqnarray}
\label{System-for-first-values-of-psi}
\psi(0)&=&\psi(1)f(0)+\overline{F}(0)\\
\nonumber
\psi(1)&=&\psi(2)f(0)+\psi(1)f(1)+\overline{F}(1)\\
\nonumber
\psi(2)&=&\psi(3)f(0)+\psi(2)f(1)+\psi(1)f(2)+\overline{F}(2)\\
\nonumber
&\vdots&\\
\nonumber
\psi(m-1)&=&\psi(m)f(0)+\psi(m-1)f(1)+\cdots+\psi(1)f(m-1)+\overline{F}(m-1).
\end{eqnarray}
This is a set of $m$ linear equations for the $m$ unknown $\psi(1),\ldots,\psi(m)$, where,
as said earlier, $\psi(0)$ is taken as known and equal to $\mu_Y$.
In matrix form, the system~\eqref{System-for-first-values-of-psi} can be
written as follows

\begin{equation}
\label{Matrix-system-for-first-values-of-psi}
\left(
\begin{array}{ccccc}
f(0) & 0 & 0 &  \cdots & 0\\
f(1)-1 & f(0) & 0 &  \cdots & 0\\
f(2) & f(1)-1 & f(0) &  \cdots & 0\\
\vdots & \vdots & \vdots & \ddots & \vdots\\
f(m-1) & f(m-2) & f(m-3) & \cdots & f(0)
\end{array}
\right)
\left(
\begin{array}{c}
\psi(1) \\
\psi(2) \\
\psi(3) \\
\vdots \\
\psi(m)
\end{array}
\right)
=
\left(
\begin{array}{c}
\psi(0)-\overline{F}(0) \\
-\overline{F}(1) \\
-\overline{F}(2) \\
\vdots \\
-\overline{F}(m-1)
\end{array}
\right).
\end{equation}
The determinant of the above matrix is $(f(0))^m$, which is strictly positive
by the net profit condition. Hence,
a unique collection of values for $\psi(1),\ldots,\psi(m)$ can be determined
from~\eqref{Matrix-system-for-first-values-of-psi}.
Given these initial values, the
recursive relation~\eqref{Recursive-relation-2} can now be used to find $\psi(u)$ for all $u\ge1$.
To this
end, define the $m$ coefficients
\begin{equation}
\label{alphas-definition}
\alpha_0=\frac{1-f(1)}{f(0)},\ \alpha_1=-\frac{f(2)}{f(0)},\ \alpha_2=-\frac{f(3)}{f(0)},\ \ldots,\ 
\alpha_{m-1}=-\frac{f(m)}{f(0)}.
\end{equation}
Observe that $\alpha_0>1$ and $\alpha_1\le0,\ldots,\alpha_{m-2}\le0$ but
$\alpha_{m-1}<0$
since $f(m)>0$.
 Hence, there is exactly one change of sign in the
sequence $\alpha_0, \alpha_1,\ldots,\alpha_{m-1}$.
Observe also that $\alpha_0+\alpha_1+\cdots+\alpha_{m-1}=1$.
Then, the recursive relation~\eqref{Recursive-relation-2} can be written as
\begin{equation}
\label{Recursive-relation-3}
\psi(u+1)=\sum_{k=0}^{m-1} \alpha_k\,\psi(u-k),\quad\mbox{for }\ u\ge m,
\end{equation}
and writing $u+m$ instead of $u$ yields
\begin{equation}
\label{Recursive-relation-4}
\psi(u+m+1)=\sum_{k=0}^{m-1} \alpha_k\,\psi(u+m-k),\quad\mbox{for }\ u\ge0.
\end{equation}
The next step is to solve the recurrence equation~\eqref{Recursive-relation-4}.
We will do that using the method of characteristic polynomials
(see~\cite{1971Brousseau, 2013Sedgewick-Flajolet}). The characteristic polynomial
associated with~\eqref{Recursive-relation-4} is
\begin{equation}
\label{Characteristic-polynomial}
p(y)=y^m-\sum_{k=0}^{m-1} \alpha_k\,y^{m-1-k}.
\end{equation}
Being $p(y)$ an $m$ degree polynomial, it has at most $m$ zeroes, which can be real or complex
and can have any multiplicity. Let us suppose that
$z_1,\ldots,z_{\ell}$ are the roots of $p(y)=0$, where $1\le\ell\le m$, with multiplicities
$n_1,\ldots,n_{\ell}$,
respectively. Observe that $n_1+\cdots+n_{\ell}=m$.
It is known~(\cite{2013Sedgewick-Flajolet}) and this is the important fact, that the general
solution to the recurrence
sequence~\eqref{Recursive-relation-4} is given by
\begin{equation}
\label{Solution-psi(u)}
\psi(u)=\sum_{k=1}^{\ell} \sum_{j=1}^{n_k} b_{k,j}\,u^{j-1}\,z_k^u,\quad\mbox{for }\ u\ge1,
\end{equation}
where $b_{k,j}$ are $m$ constants which can be chosen so that the
expression~\eqref{Solution-psi(u)} matches
the initial data $\psi(1),\ldots,\psi(m)$. From~\eqref{Solution-psi(u)}, observe that
each root $z_k$ with multiplicity $n_k$ will have $n_k$ associated
constants $b_{k,1},\ldots,b_{k,n_k}$. When the multiplicity of $z_k$  is $1$ there is only one
constant and we write $b_k$ instead of $b_{k,1}$.
For example, suppose that $m=6$ and that we have $2$ simple positive roots $z_1$ and $z_2$,
and $2$ roots $z_3$ and $z_4$ each with multiplicity $2$. The $6$ constants $b_{k,j}$
of~\eqref{Solution-psi(u)} that satisfy the initial conditions are given by the solution to the
system	

\begin{equation}
\label{Example-1-matrix-Z-system}
\left(
\begin{array}{c|c|cc|cc}
z_1   & z_2   & z_3    & z_3    & z_4    & z_4   \\
z_1^2 & z_2^2 & z_3^2 & 2z_3^2 & z_4^2 & 2z_4^2 \\
z_1^3 & z_2^3 & z_3^3 & 3z_3^3 & z_4^3& 3z_4^3 \\
z_1^4 & z_2^4 & z_3^4 & 4z_3^4& z_4^4& 4z_4^4 \\
z_1^5 & z_2^5 & z_3^5 & 5z_3^5 & z_4^5& 5z_4^5 \\
z_1^6 & z_2^6 & z_3^6 & 6z_3^6 & z_4^6& 6z_4^6 \\
\end{array}
\right)\left(
\begin{array}{c}
b_1\\
\mbox{-----}\\
b_{2}\\
\mbox{-----}\\
b_{3,1}\\
b_{3,2}\\
\mbox{-----}\\
b_{4,1}\\
b_{4,2}
\end{array}
\right)=\left(
\begin{array}{c}
\psi(1)\\
\psi(2)\\
\psi(3)\\
\psi(4)\\
\psi(5)\\
\psi(6)\\
\end{array}
\right).
\end{equation}
The vertical and horizontal lines help visually separate the terms associated with
the different roots. 
In the general case, the $m$ constants $b_{k,j}$ are the solution to the
system of linear equations
\begin{equation}
\label{system2-for-b}
\mathbf{Z}
\left(
\begin{array}{l}
b_{1}\\
b_{2}\\
b_{3,-}\\
\vdots\\
b_{\ell,-}
\end{array}
\right)
=
\left(
\begin{array}{c}
\psi(1)\\
\psi(2)\\
\psi(3)\\
\vdots\\
\psi(m)
\end{array}
\right),
\end{equation}
with $b_{k,-}$ being the column vector $(b_{k,1},\ldots,b_{k,n_k})^t$ for $k=3,\ldots,\ell$,
and $\mathbf{Z}$ is the $m\times m$ matrix of the form
$(\mathbf{\underline{Z}_1},\mathbf{\underline{Z}_2},\ldots,\mathbf{\underline{Z}_\ell})$,
where the $k$-entry is given by the following $m\times n_k$ submatrix
\begin{equation*}
\quad
\mathbf{\underline{Z}}_k=\left(
\begin{array}{cccc}
z_k & z_k & \cdots & z_k \\
z_k^2 & 2z_k^2 & \cdots & 2^{n_k-1}z_k^2 \\
z_k^3 & 3z_k^3 & \cdots & 3^{n_k-1}z_k^3 \\
\vdots & \vdots & \vdots & \vdots \\
z_k^m & mz_k^m & \cdots & m^{n_k-1}z_k^m \\
\end{array}
\right)_{m\times n_k}.
\end{equation*}

Observe that each root $z_k$ has the associated submatrix $\mathbf{\underline{Z}}_k$.
It is clear that when a root $z_k$ is simple ($n_k=1$), the submatrix
$\mathbf{\underline{Z}}_k$ is a column vector.
As a summary, we now write the statement of our main result.

\begin{theorem}
\label{Main-theorem}
For the Gerber-Dickson risk process~\eqref{Gerber-Dickson-model} with claims $Y$
having a probability function $f(x)$ with support $\{0,1,\ldots,m\}$ with $m\ge2$ and
such that $E(Y)<1$, the ultimate ruin probability $\psi(u)$, for $u\ge1$, satisfies the recurrence
sequence~\eqref{Recursive-relation-4} and can be written as in~\eqref{Solution-psi(u)}.
\end{theorem}

In a nutshell, the general procedure to obtain the ruin probability values in~\eqref{Solution-psi(u)} is as follows:
[a] Given a probability function $f(x)$ for the claims, calculate the coefficients
$\alpha_0, \alpha_1,\ldots, \alpha_{m-1}$ using~\eqref{alphas-definition} and construct the
characteristic polynomial~\eqref{Characteristic-polynomial}.
[b] Find the roots $z_1,\ldots,z_m$ of~\eqref{Characteristic-polynomial}.
[c] Compute the initial values $\psi(1),\ldots,\psi(m)$ using~\eqref{Matrix-system-for-first-values-of-psi}.
[d] Finally, find the constant $b_{k,j}$ using~\eqref{system2-for-b}.
These  steps yield all the terms appearing in our general formula~\eqref{Solution-psi(u)}.\\

As a special case, suppose that all the roots of the polynomial~\eqref{Characteristic-polynomial}
are simple, then
there are $m$ different roots $z_1,\ldots,z_m$ each with multiplicity $1$ and the ruin probability
formula~\eqref{Solution-psi(u)} reduces to
\begin{equation}
\label{Solution-for-simple-roots}
\psi(u)=\sum_{k=1}^{m} b_{k}\,z_k^u,\quad\mbox{for }\ u\ge1,
\end{equation}
where $b_1,\ldots,b_m$ are the solution to the linear system
\begin{equation}
\label{system-for-b-simple-case}
\left(
\begin{array}{cccc}
z_1 & z_2 & \cdots & z_m\\
z_1^2 & z_2^2 & \cdots & z_m^2\\
\vdots & \vdots & \ddots & \vdots\\
z_1^{m} & z_2^{m} & \cdots & z_m^{m}
\end{array}
\right)
\left(
\begin{array}{c}
b_1 \\
b_2 \\
\vdots \\
b_m
\end{array}
\right)
=
\left(
\begin{array}{c}
\psi(1) \\
\psi(2) \\
\vdots \\
\psi(m)
\end{array}
\right).
\end{equation}

In this way we have translated a general problem from the theory of ruin into two classical mathematical
problems: finding the roots of a polynomial and solving a system of linear equations. Of course, the
main difficulty in applying this method lies in finding the roots of the polynomial.

\begin{example}
\label{Player-example} (A gambler's ruin problem)
Assume that the distribution of the claims $Y$ is given by
\begin{displaymath}
f(y)=\left\{
\begin{array}{ll}
p & \mbox{if }\ y=0,\\
q &\mbox{if }\ y=2,\\
0 & \mbox{otherwise,}
\end{array}
\right.
\end{displaymath}
where $p+q=1$, that is, $m=2$ and the claims can take only two values: $0$ or $2$.
The net profit condition reads
$E(Y)=2(1-p)<1$, i.e. $p>1/2$. Since the insurance premium is $1$ in each period, the reserve
process~\eqref{Gerber-Dickson-model} goes up $1$ unit with probability $p$ or goes down $1$
unit with probability $q$ at the end of each unit time period. This is the same situation of a
player $A$ with initial capital $u\ge1$ who sequentially bets $1$ unit of currency against player $B$,
where on each trial his probability of winning is $p$ and his probability of losing is $q$.
Assuming that player $B$ has an infinite amount of money, it can be
proved~(\cite[p.~143]{taylor1998introduction})
that the ultimate ruin probability for player $A$ is
\begin{equation}
\label{Player-A-example-ultimate-ruin-prob}
P(``\mbox{Ultimate ruin}")=({q}/{p})^u,\quad\mbox{for }\ u\ge1.
\end{equation}
\begin{enumerate}[a)]
\item We can recover formula~\eqref{Player-A-example-ultimate-ruin-prob} using Theorem~\ref{Main-theorem} as follows. We have $m=2$ and
the constants $\alpha_k$
defined in~\eqref{alphas-definition} are given by $\alpha_0=1/p$ and $\alpha_1=-q/p$.
The characteristic polynomial~\eqref{Characteristic-polynomial} is
$p(y)=y^2-(1/p)y+q/p=(y-1)(y-q/p)$
with two different real roots $z_1=1$ and $z_2=q/p$. The initial values $\psi(1)$ and $\psi(2)$ are
given by the system~\eqref{Matrix-system-for-first-values-of-psi}, which in this case is
\begin{equation}
\left(
\begin{array}{cc}
p & 0 \\
-1 & p \\
\end{array}
\right)
\left(
\begin{array}{c}
\psi(1) \\
\psi(2) \\
\end{array}
\right)
=
\left(
\begin{array}{c}
q \\
-q \\
\end{array}
\right),
\end{equation}
with solution $\psi(1)=q/p$ and $\psi(2)=(q/p)^2$. 
The general ruin probability formula~\eqref{Solution-for-simple-roots}
is $\psi(u)=b_1z_1^u+b_2z_2^u=b_1+b_2(q/p)^u$, for $u\ge1$,
where the coefficients $b_1$ and $b_2$ are such that the initial values
$\psi(1)=q/p$ and $\psi(2)=(q/p)^2$ are satisfied. This leads to $b_1=0$ and $b_2=1$ and
the formula~\eqref{Player-A-example-ultimate-ruin-prob} is thus rediscovered.

\item Generating functions can also be used to solve the gambler's ruin problem and find again
formula~\eqref{Player-A-example-ultimate-ruin-prob}. The recurrence sequence for this problem is
given by
\begin{equation}
\label{Recurrence-eqn-for-players-ruin-problem}
\psi(u)=\psi(u+1)f(0)+\psi(u-1)f(2),\quad\mbox{for }\ u\ge1,
\end{equation}
with initial values $\psi(1)=q/p$ and $\psi(2)=(q/p)^2$. Define $G(s)=\sum_{u=1}^\infty \psi(u)\,s^u$.
Multiplying~\eqref{Recurrence-eqn-for-players-ruin-problem} by $s^u$ and summing gives
\begin{equation}
\label{Eqn-for-players-ruin-problem}
\sum_{u=2}^\infty \psi(u)\,s^u=f(0)\sum_{u=2}^\infty \psi(u+1)\,s^u+f(2)\sum_{u=2}^\infty \psi(u-1)\,s^u.
\end{equation}
From~\eqref{Eqn-for-players-ruin-problem} an equation for $G(s)$ can then be derived and after
some simplifications using the initial data  the result below is obtained. The coefficients in the last sum
are the solution $\psi(u)$, for $u\ge1$.
\begin{displaymath}
G(s)=\frac{(q/p)s}{1-(q/p)s}=\sum_{u=1}^\infty (q/p)^us^u.
\end{displaymath}
\end{enumerate}
\end{example}

\section{Discussion on the main result}
In this section we  make some comments on the main formula~\eqref{Solution-psi(u)}
and its components. Despite the fact that some of the roots $z_k$ and the coefficients $b_{k,j}$
in~\eqref{Solution-psi(u)} can be complex, the whole expression for $\psi(u)$ always gives a real number in $(0,1)$, as expected. We will explain why this is so below.

\subsection{On the roots of the characteristic polynomial}

It is clear that finding the roots of the characteristic polynomial is the main difficulty in
applying formula~\eqref{Solution-psi(u)}. In this section we apply some basic
results on polynomials to shed some light on this problem.
In the following section we will show numerical cases where the roots are found using a computer.

\begin{enumerate}
\item Since $\alpha_0+\alpha_1+\cdots+\alpha_{m-1}=1$, $z_1=1$ is always a root
of~\eqref{Characteristic-polynomial}. Observe that the coefficient $b_1$ associated with $z_1=1$
must be zero since the contribution of $z_1$ to the ruin probability reduces to $\psi_1(u)=b_1$, which
should converge to zero as $u\to\infty$. This is only possible when $b_1=0$.
%
%
\item The sequence of coefficients $1,-\alpha_0,-\alpha_1,\ldots,-\alpha_{m-2},-\alpha_{m-1}$ of the
characteristic polynomial~\eqref{Characteristic-polynomial} has exactly $2$ changes of sign. By the
Descartes' rule of signs~(\cite{doi:10.1080/00029890.1998.12004907}), the number of positive
roots of~\eqref{Characteristic-polynomial}, counted with multiplicity, is $0$ or $2$. 
Since $z_1=1$ is always a root, we conclude that there are exactly $2$ positive roots which are
denoted by $z_1=1$ and $z_2>0$.
\item Since $z_1=1$ is always a root, simple calculations show that $p(y)=(y-1)q(y)$ where
\begin{equation}
\label{Polynomial-q}
q(y)=y^{m-1}-\sum_{k=1}^{m-1}\,\frac{\overline{F}(y)}{f(0)}\,y^{m-k-1}.
\end{equation}
The coefficients $\overline{F}(y)/f(0)$ are non-negative and at least one of them is nonzero.
In fact, they are all strictly positive since $f(m)>0$.
By the Cauchy's theorem on the roots of polynomials~(\cite{2010Prasolov}), $q(y)$ has a unique simple positive
root $z_2$ and the norm of any other root $z$ is such that $|z|\le z_2$. We already knew the positivity
of $z_2$ from Descartes' rule of signs.
\item One the many bounding results known in the literature~(\cite{Hirst-Macey-1997}) on the roots of
polynomials implies that any root $z$ of $q(y)$ satisfies
\begin{equation}
\label{Upper-bound-1-for roots}
|z|\le\max\,\{\,1,\,\sum_{k=1}^{m-1}\frac{\overline{F}(y)}{f(0)}\,\}.
\end{equation}
Inasmuch as
\begin{displaymath}
\sum_{k=1}^m\frac{\overline{F}(y)}{f(0)}=\frac{1}{f(0)}[E(Y)-\overline{F}(0)]=\frac{1}{f(0)}[E(Y)-(1-f(0))]
=1-\frac{1-E(Y)}{f(0)}<1,
\end{displaymath}
we have that the right-hand side of~\eqref{Upper-bound-1-for roots} is $1$ and thus, 
$z_2\le1$. Since $z_2$ cannot be equal to $z_1=1$, we conclude that any root $z$ (different from $z_1$) of the
characteristic polynomial $p(y)$ is such that 
\begin{equation}
|z|\le z_2<z_1=1.
\end{equation}
This means all roots different from $z_1$ lie inside the circle of radius $z_2$ in the complex plane.
\item Since the coefficients of the polynomial are all real, roots occur in conjugate pairs,
i.e. if $z$ is a root, then its complex conjugate $\overline{z}$ is also a root. This is known as the
complex conjugate root theorem. Furthermore, it can be
shown that two conjugate roots $z_k$ and $\overline{z}_k$ share the same multiplicity $n_k$.
Moreover, it can also be proved that the list of coefficients $b_{k,j}$ associated with two conjugate roots
$z_k$ and $\overline{z}_k$ are also conjugate. All these results imply that the formula~\eqref{Solution-psi(u)}
yields a real number for $\psi(u)$.
\item From the above results it follows that when $m\ge2$ is odd there is always a
negative root.
\end{enumerate}

\subsection{An alternative recurrence sequence}
The second recurrence equation~\eqref{Recursive-equation-2} for $\psi(u)$ can also be used as
the starting point in our procedure. We here show that both schemes, 
\eqref{Recursive-equation-1} and~\eqref{Recursive-equation-2}, are the same as expected.
For $u\ge m$, the second sum on the right-hand side
of~\eqref{Recursive-equation-2} vanishes and we have the relation
\begin{displaymath}
\psi(u)=\sum\limits_{k=0}^{m-1}\,\overline{F}(k)\,\psi(u-k),
\quad\mbox{for } u\ge m,
\end{displaymath}
where the upper limit of the sum is replaced by $m-1$. Since $1-\overline{F}(0)=f(0)$, solving for
$\psi(u)$ and writing $u+m$ instead of $u$ yields
\begin{displaymath}
\psi(u+m)=\sum\limits_{k=1}^{m-1}\,\frac{\overline{F}(k)}{f(0)}\,\psi(u+m-k),
\quad\mbox{for } u\ge0.
\end{displaymath}
This is a recurrence sequence of order $m-1$ with characteristic polynomial $q(y)$ defined
in~\eqref{Polynomial-q} and the relation  $(y-1)q(y)=p(y)$ can be verified.
Thus, leaving aside the root $z_1=1$, the
polynomials $q(y)$ and $p(y)$ have exactly the same $m-1$ roots and the same
formula~\eqref{Solution-psi(u)} for the ruin probability is again obtained.
Since the recursive equations~\eqref{Recursive-equation-1} and~\eqref{Recursive-equation-2} are
equivalent, the initial values $\psi(1),\ldots,\psi(m)$ are the same and the linear
system~\eqref{Matrix-system-for-first-values-of-psi} produces the same associated coefficients $b_{k,j}$.

\section{An approximation}
The numerical examples presented in the next section show that the second positive root $z_2$
and its associated coefficient $b_2$ play a leading role in the ruin probability formula~\eqref{Solution-psi(u)}.
This justifies looking for an approximate solution to the linear system~\eqref{system2-for-b}  of the
form $(0,b_2,0,\ldots,0)$, which gives the following approximation for the ruin probability.

\begin{corollary}
For the Gerber-Dickson risk process~\eqref{Gerber-Dickson-model},
\begin{equation}
\label{Approximation-1}
\psi(u)\approx\hat{\psi}_1(u):=b_2\,z_2^u,\quad\mbox{for }\ u\ge1.
\end{equation}
\end{corollary}

This simple approximation is the contribution of the second positive root $z_2$ to the ruin probability
formula~\eqref{Solution-psi(u)}. It only works when one has some knowledge on the values of $b_2$ and $z_2$.
We here assume $b_2>0$. When $z_2$ and $b_2$ are completely unknown (as it is in most cases), they can be
approximated as follows.
The first two equations of the system~\eqref{system2-for-b} for an
approximate solution of the form  $(0,b_2,0,\ldots,0)$ are
\begin{eqnarray*}
b_2\,z_2&=&\psi(1),\\
b_2\,z_2^2&=&\psi(2),
\end{eqnarray*}
which yields $b_2=(\psi(1))^2/\psi(2)$ and $z_2=\psi(2)/\psi(1)$.
Substituting these values in~\eqref{Approximation-1} gives the following estimate.
\begin{corollary}
For the Gerber-Dickson risk process~\eqref{Gerber-Dickson-model},
\begin{equation}
\label{Approximation-2}
\psi(u)\approx\hat{\psi}_2(u):=
\psi(1)\left(\frac{\psi(2)}{\psi(1)}\right)^{u-1},\quad\mbox{for }\ u\ge1.
\end{equation}
\end{corollary}
Observe that the exact values of $\psi(1)$ and $\psi(2)$
can be obtained from~\eqref{System-for-first-values-of-psi}, namely,
\begin{eqnarray}
\label{psi(1)-formula}
\psi(1)&=&1-\frac{1-E(Y)}{f(0)},\\
\label{psi(2)-formula}
\psi(2)&=&1-\frac{1-E(Y)}{f(0)}\cdot\frac{1-f(1)}{f(0)}.
\end{eqnarray}
It is straightforward to check from~\eqref{Approximation-2} that $\psi(1)=\hat{\psi}_2(1)$ and
$\psi(2)=\hat{\psi}_2(2)$, i.e. the approximation $\hat{\psi}_2(u)$ is exact for $u=1,2$ and is
based on only three quantities from the claims distribution: $f(0)$, $f(1)$ and $E(Y)$.\\
\begin{example}
When claims follow the geometric distribution $f(y)=p(1-p)^y$ for $y=0,1,\ldots$, the
net profit condition reads $p>1/2$. The initial ruin probability values are $\psi(1)=((1-p)/p)^2$ and
$\psi(2)=((1-p)/p)^3$. Although this distribution fails to have bounded support, the
approximation~\eqref{Approximation-2} yields the exact ruin probability~(see \cite{2017dickson}),
\begin{equation}
\label{psi(u)-geo-dist}
\psi(u)=\hat{\psi}_2(u)=\left(\frac{1-p}{p}\right)^{u+1},\quad\mbox{for }\ u\ge0.
\end{equation}
\end{example}
\begin{example}
Assume that the claims $Y$ follows a distribution within the class $(a,b,0)$, that is, the relation
$f(y)=(a+b/y)f(y-1)$ is satisfied for $y=1,2,\ldots$, for real constants $a$ and $b$, see~\cite{2017dickson}.
A distribution within this class has support $\{0,1,2,\ldots\}$ or $\{0,1,\ldots,n\}$ for some
integer $n$.
It is well known that this class of distributions only includes the binomial, the negative binomial and
the Poisson distribution for certain values of $a$ and $b$. The geometric distribution is a particular case of the negative binomial distribution and it was  mentioned in the previous example.
It can be shown that
$f(0)=(1-a)^{(a+b)/a}$, $f(1)=(a+b)(1-a)^{(a+b)/a}$ and $E(Y)=(a+b)/(1-a)$ for $a\ne1$.
The net profit condition requires that the values of  $a$ and $b$ must be such that $(a+b)/(1-a)<1$.
With this information and using~\eqref{psi(1)-formula} and~\eqref{psi(2)-formula}, the values
for $\psi(1)$ and $\psi(2)$ can be calculated. Formula~\eqref{Approximation-2}
then yields the following approximate ruin probability formula
\begin{equation}
\label{Approximation-2-for-(a,b,0)}
\hat{\psi}_2(u)=
\left(1-\frac{1-(a+b)/(1-a)}{(1-a)^{(a+b)/a}}\right)\left(\frac{1-\frac{1-(a+b)/(1-a)}{(1-a)^{(a+b)/a}}\cdot\frac{1-(a+b)(1-a)^{(a+b)/a}}{(1-a)^{(a+b)/a}}}{1-\frac{1-(a+b)/(1-a)}{(1-a)^{(a+b)/a}}}\right)^{u-1},
\end{equation}
for $u\ge1$ and for any claims distribution within the $(a,b,0)$ class.
Formula~\eqref{Approximation-2-for-(a,b,0)} can be further reduced according to the particular
values of $a$ and $b$. For example, when $a=1-p$ and $b=0$ the geometric distribution
is obtained and~\eqref{Approximation-2-for-(a,b,0)} reduces to~\eqref{psi(u)-geo-dist}.
When $a\to0$ and for $b=\lambda>0$, the $Poisson(\lambda)$ distribution is obtained and a
shorter formula can be derived.
\end{example}
In the following section we will numerically show how good the approximations $\hat{\psi}_1(u)$
and $\hat{\psi}_2(u)$ are in some particular cases.

\section{Numerical examples}

We now present some numerical examples of particular distributions for the claims
which yield exact ruin probabilities $\psi(u)$ using formula~\eqref{Solution-psi(u)}.
The roots of the characteristic polynomial~\eqref{Characteristic-polynomial} and the solution
to the linear system~\eqref{system2-for-b} are found using the
{\tt R} pracma package~(\cite{R-pracma}).
In each case observe that the approximation $\hat{\psi}_1(u)=b_2\,z_2^u$,
for $u\ge1$, seems to be accurate. This means that the term $b_2z_2^u$ is the one making the
major contribution to the exact ruin probability formula~\eqref{Solution-psi(u)} for $\psi(u)$.\\

It is apparent that a correct calculation of the roots of~\eqref{Characteristic-polynomial} is absolutely essential for providing exact ruin probabilities
trough formula~\eqref{Solution-psi(u)}. When using a computer, numerical errors can lead to inaccurate ruin probabilities. For example, in our numerical experimentation we found that for claims with
small mean $\mu_Y$, the function {\tt polyroots} of the {\tt R} pracma package,
sometimes produced some lack of  precision in the
calculations of the roots of the characteristic polynomial. Hence, for any numerical application
of the ruin probability formula~\eqref{Solution-psi(u)}, it is advisable to make sure the roots
are properly calculated.
In the examples shown below and to save some space, most of the numbers are
rounded up to three or four decimal digits.\\



\example{
\label{Numerical-example-1}
Suppose that $m=2$ and claims $Y$ follow the distribution $f(0)={1}/{2}$, $f(1)={1}/{4}$ and $f(2)={1}/{4}$, with mean $\mu_Y=\psi(0)=3/4$. The coefficients $\alpha$ can then be calculated
and the characteristic polynomial can be constructed.
The two roots and their multiplicities are shown in Figure~\ref{Figure1-example1}.\vspace*{-1.3cm}

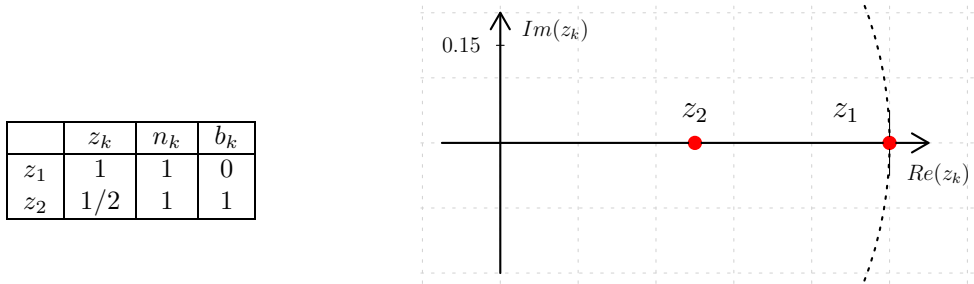
\begin{figure}[H]
\begin{floatrow}
\capbtabbox{%
\footnotesize
%
\begin{tabular}{|c|c|c|c|}
\hline
 & $z_k$ & $n_k$ &  $b_k$\\ 
\hline
$z_1$ & $1$ & $1$ & $0$ \\
$z_2$ & $1/2$ & $1$ & $1$\\
\hline
\end{tabular}
\vspace{3cm}
}{%
}
\begin{tikzpicture}[x=1pt,y=1pt]
\definecolor{fillColor}{RGB}{255,255,255}
\path[use as bounding box,fill=fillColor,fill opacity=0.00] (0,0) rectangle (289.08,216.81);
\begin{scope}
\path[clip] ( 49.20, 61.20) rectangle (263.88,167.61);
\definecolor{fillColor}{RGB}{255,0,0}

\path[fill=fillColor] (152.86,114.41) circle (  2.70);

\path[fill=fillColor] (226.48,114.41) circle (  2.70);
\end{scope}
\begin{scope}
\path[clip] ( 49.20, 61.20) rectangle (263.88,167.61);
\definecolor{drawColor}{RGB}{211,211,211}

\path[draw=drawColor,line width= 0.4pt,dash pattern=on 1pt off 3pt ,line join=round,line cap=round] ( 49.79, 61.20) -- ( 49.79,167.61);

\path[draw=drawColor,line width= 0.4pt,dash pattern=on 1pt off 3pt ,line join=round,line cap=round] ( 79.24, 61.20) -- ( 79.24,167.61);

\path[draw=drawColor,line width= 0.4pt,dash pattern=on 1pt off 3pt ,line join=round,line cap=round] (108.69, 61.20) -- (108.69,167.61);

\path[draw=drawColor,line width= 0.4pt,dash pattern=on 1pt off 3pt ,line join=round,line cap=round] (138.13, 61.20) -- (138.13,167.61);

\path[draw=drawColor,line width= 0.4pt,dash pattern=on 1pt off 3pt ,line join=round,line cap=round] (167.58, 61.20) -- (167.58,167.61);

\path[draw=drawColor,line width= 0.4pt,dash pattern=on 1pt off 3pt ,line join=round,line cap=round] (197.03, 61.20) -- (197.03,167.61);

\path[draw=drawColor,line width= 0.4pt,dash pattern=on 1pt off 3pt ,line join=round,line cap=round] (226.48, 61.20) -- (226.48,167.61);

\path[draw=drawColor,line width= 0.4pt,dash pattern=on 1pt off 3pt ,line join=round,line cap=round] (255.93, 61.20) -- (255.93,167.61);

\path[draw=drawColor,line width= 0.4pt,dash pattern=on 1pt off 3pt ,line join=round,line cap=round] ( 49.20, 65.14) -- (263.88, 65.14);

\path[draw=drawColor,line width= 0.4pt,dash pattern=on 1pt off 3pt ,line join=round,line cap=round] ( 49.20, 89.77) -- (263.88, 89.77);

\path[draw=drawColor,line width= 0.4pt,dash pattern=on 1pt off 3pt ,line join=round,line cap=round] ( 49.20,114.41) -- (263.88,114.41);

\path[draw=drawColor,line width= 0.4pt,dash pattern=on 1pt off 3pt ,line join=round,line cap=round] ( 49.20,139.04) -- (263.88,139.04);

\path[draw=drawColor,line width= 0.4pt,dash pattern=on 1pt off 3pt ,line join=round,line cap=round] ( 49.20,163.67) -- (263.88,163.67);
\definecolor{drawColor}{RGB}{0,0,0}

\path[draw=drawColor,line width= 0.8pt,line join=round,line cap=round] ( 57.15,114.41) -- (241.20,114.41);

\path[draw=drawColor,line width= 0.8pt,line join=round,line cap=round] (234.95,110.79) --
	(241.20,114.41) --
	(234.95,118.02);

\path[draw=drawColor,line width= 0.8pt,line join=round,line cap=round] ( 79.24, 65.14) -- ( 79.24,163.67);

\path[draw=drawColor,line width= 0.8pt,line join=round,line cap=round] ( 82.85,157.41) --
	( 79.24,163.67) --
	( 75.62,157.41);

\path[draw=drawColor,line width= 0.8pt,dash pattern=on 1pt off 3pt ,line join=round,line cap=round] ( 79.24,114.41) circle (147.24);

\node[text=drawColor,anchor=base,inner sep=0pt, outer sep=0pt, scale=  0.70] at (245,100.34) {$Re(z_k)$};

\node[text=drawColor,anchor=base,inner sep=0pt, outer sep=0pt, scale=  0.70] at ( 100,155) {$Im(z_k)$};

\path[draw=drawColor,line width= 0.4pt,line join=round,line cap=round] (226.48,102.09) --
	(226.48,126.72);

\path[draw=drawColor,line width= 0.4pt,line join=round,line cap=round] ( 77.77,151.35) --
	( 80.71,151.35);

\node[text=drawColor,anchor=base,inner sep=0pt, outer sep=0pt, scale=  0.70] at ( 64.51,148.92) {$0.15$};
\definecolor{fillColor}{RGB}{255,0,0}

\path[fill=fillColor] (152.86,114.41) circle (  2.25);

\path[fill=fillColor] (226.48,114.41) circle (  2.25);

\node[text=drawColor,anchor=base,inner sep=0pt, outer sep=0pt, scale=  1.00] at (152.86,124.14) {$z_2$};

\node[text=drawColor,anchor=base,inner sep=0pt, outer sep=0pt, scale=  1.00] at (210,124.14) {$z_1$};
\end{scope}
\end{tikzpicture}
\end{floatrow}
\vspace{-1.5cm}
\caption{\label{Figure1-example1}
[Example~\ref{Numerical-example-1}]
Zeroes of  the characteristic polynomial.}
\end{figure}

Knowing the values of the roots $z_1=1$ and $z_2=1/2$,
the linear system~\eqref{system2-for-b} can now be numerically
solved to obtain the associated coefficients $b_1=0$ and $b_2=1$.
Since $z_2$ is the only root different from $z_1$,  it can be easily shown that
$\psi(u)=\hat{\psi}_1(u)=\hat{\psi}_1(u)=(1/2)^u$ for $u\ge1$.
These ruin probabilities are shown in Figure~\ref{Figure2-example1}.
\vspace*{-1.3cm}

\begin{figure}[H]
\begin{floatrow}
\capbtabbox{%
\footnotesize
%
\begin{tabular}{|c|l|l|l|}
\hline
$u$ & $\psi(u)$ & $\hat{\psi}_1(u)$ & $\hat{\psi}_2(u)$ \\ 
\hline
$0$ & $0.75$ & $-$ & $-$ \\
$1$ & $0.5$ & $0.5$ & $0.5$ \\
$2$ & $0.25$ & $0.25$ & $0.25$ \\
$3$ & $0.125$ & $0.125$ & $0.125$ \\
$4$ & $0.0625$ & $0.0625$ & $0.0625$ \\
$5$ & $0.03125$ & $0.03125$ & $0.03125$ \\
$6$ & $0.015625$ & $0.015625$ & $0.015625$ \\
\hline
\end{tabular}
\vspace{2cm}
}{%
}
\begin{tikzpicture}[x=1pt,y=1pt]
\definecolor{fillColor}{RGB}{255,255,255}
\path[use as bounding box,fill=fillColor,fill opacity=0.00] (0,0) rectangle (289.08,216.81);
\begin{scope}
\path[clip] ( 49.20, 61.20) rectangle (263.88,167.61);
\definecolor{drawColor}{RGB}{0,0,0}

\path[draw=drawColor,line width= 1.6pt,line join=round,line cap=round] ( 95.01,158.98) --
	(132.88,112.06) --
	(170.74, 88.60) --
	(208.60, 76.87) --
	(246.46, 71.01);
\definecolor{drawColor}{RGB}{255,0,0}

\path[draw=drawColor,line width= 1.2pt,dash pattern=on 4pt off 4pt ,line join=round,line cap=round] ( 95.01,158.98) --
	(132.88,112.06) --
	(170.74, 88.60) --
	(208.60, 76.87) --
	(246.46, 71.01);
\definecolor{drawColor}{RGB}{0,0,255}

\path[draw=drawColor,line width= 1.2pt,dash pattern=on 1pt off 3pt ,line join=round,line cap=round] ( 95.01,158.98) --
	(132.88,112.06) --
	(170.74, 88.60) --
	(208.60, 76.87) --
	(246.46, 71.01);
\definecolor{drawColor}{RGB}{0,0,0}

\path[draw=drawColor,line width= 0.3pt,line join=round,line cap=round] ( 57.15, 65.14) -- (255.93, 65.14);

\path[draw=drawColor,line width= 0.3pt,line join=round,line cap=round] (249.67, 61.53) --
	(255.93, 65.14) --
	(249.67, 68.75);

\path[draw=drawColor,line width= 0.3pt,line join=round,line cap=round] ( 57.15, 65.14) -- ( 57.15,163.67);

\path[draw=drawColor,line width= 0.3pt,line join=round,line cap=round] ( 60.76,157.41) --
	( 57.15,163.67) --
	( 53.54,157.41);

\node[text=drawColor,anchor=base,inner sep=0pt, outer sep=0pt, scale=  1.00] at (190,150) {$\psi=\hat{\psi}_1=\hat{\psi}_2$};

\node[text=drawColor,anchor=base,inner sep=0pt, outer sep=0pt, scale=  1.00] at (163,135) {$\hat{\psi}_1$};
\node[text=drawColor,anchor=base,inner sep=0pt, outer sep=0pt, scale=  1.00] at (163,120) {$\hat{\psi}_2$};

\definecolor{drawColor}{RGB}{0,0,255}
\path[draw=drawColor,line width= 1.2pt,dash pattern=on 1pt off 3pt ,line join=round,line cap=round] (175,138) --
	(197,138);

\definecolor{drawColor}{RGB}{225,0,0}
\path[draw=drawColor,line width= 1.2pt,dash pattern=on 4pt off 4pt ,line join=round,line cap=round] (175,125) --
	(197,125);

\definecolor{drawColor}{RGB}{0,0,0}
\node[text=drawColor,anchor=base,inner sep=0pt, outer sep=0pt, scale=  1.00] at (255.93, 71.05) {$u$};
\end{scope}
\begin{scope}
\path[clip] (  0.00,  0.00) rectangle (289.08,216.81);
\definecolor{drawColor}{RGB}{0,0,0}

\node[text=drawColor,anchor=base,inner sep=0pt, outer sep=0pt, scale=  1.00] at ( 57.15, 51.60) {0};

\node[text=drawColor,anchor=base,inner sep=0pt, outer sep=0pt, scale=  1.00] at ( 95.01, 51.60) {1};

\node[text=drawColor,anchor=base,inner sep=0pt, outer sep=0pt, scale=  1.00] at (132.88, 51.60) {2};

\node[text=drawColor,anchor=base,inner sep=0pt, outer sep=0pt, scale=  1.00] at (170.74, 51.60) {3};

\node[text=drawColor,anchor=base,inner sep=0pt, outer sep=0pt, scale=  1.00] at (208.60, 51.60) {4};

\node[text=drawColor,anchor=base,inner sep=0pt, outer sep=0pt, scale=  1.00] at (246.46, 51.60) {5};

\node[text=drawColor,anchor=base east,inner sep=0pt, outer sep=0pt, scale=  1.00] at ( 55.20, 99.23) {0.2};

\node[text=drawColor,anchor=base east,inner sep=0pt, outer sep=0pt, scale=  1.00] at ( 55.20,136.77) {0.4};
\end{scope}
\end{tikzpicture}
\end{floatrow}
\vspace{-1.5cm}
\caption{\label{Figure2-example1}
[Example~\ref{Numerical-example-1}]
Ruin probabilities and their approximations.}
\end{figure}
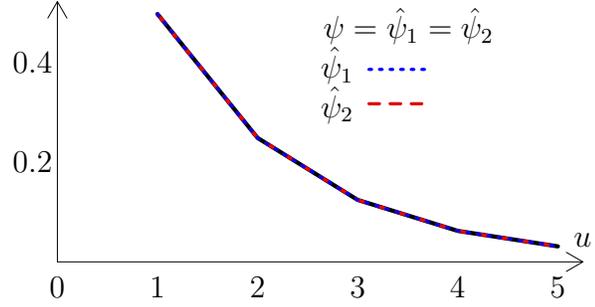
} 



\example{
\label{Numerical-example-2}
Suppose that $m=5$ and claims $Y$ follow a Binomial$(5,99/500)$ distribution with mean
$\mu_Y=\psi(0)=0.99$. The coefficients $\alpha$ are calculated and the characteristic
polynomial~\eqref{Characteristic-polynomial} can be constructed. The five
roots and their multiplicities are shown in Figure~\ref{Figure1-example2}.
We have two positive roots $z_1$ and $z_2$, one negative root $z_3$, and two complex roots $z_4$ and $z_5$,
which are complex conjugate. All roots have multiplicity $1$.\vspace*{-0.6cm}

\begin{figure}[H]
\begin{floatrow}
\capbtabbox{%
\footnotesize
%
\begin{tabular}{|c|l|l|l|}
\hline
 & $z_k$ & $n_k$ & $b_{k,j}$\\ 
\hline
$z_1$ & $1$ & $1$ & $0$ \\
$z_2$ & $0.975$ & $1$ & $0.995$\\
$z_3$ & $-0.080$ & $1$ & $1.556E-03$\\
$z_4$ & $-0.057+0.091i$ & $1$ & $(1.721+1.025i)E-03$\\
$z_5$ & $-0.057-0.091i$ & $1$ & $(1.721-1.025i)E-03$\\
\hline
\end{tabular}
\vspace{1.5cm}
}{%
}
\begin{tikzpicture}[x=0.5pt,y=0.6pt, base/.style = {rectangle, draw, align=center}]
\definecolor{fillColor}{RGB}{255,255,255}
\path[use as bounding box,fill=fillColor,fill opacity=0.00] (0,0) rectangle (361.35,289.08);
\begin{scope}
\path[clip] (0,0) rectangle (336.15,239.88);
\definecolor{fillColor}{RGB}{255,0,0}
\path[fill=fillColor] ( 73.49,150.54) circle (  2.70);

\path[fill=fillColor] ( 78.02,112.74) circle (  2.70);

\path[fill=fillColor] ( 78.02,188.34) circle (  2.70);

\path[fill=fillColor] (281.29,150.54) circle (  2.70);

\path[fill=fillColor] (286.16,150.54) circle (  2.70);
\end{scope}
\begin{scope}
\path[clip] ( 49.20, 61.20) rectangle (336.15,239.88);
\definecolor{drawColor}{RGB}{211,211,211}

\path[draw=drawColor,line width= 0.4pt,dash pattern=on 1pt off 3pt ,line join=round,line cap=round] ( 49.99, 61.20) -- ( 49.99,239.88);

\path[draw=drawColor,line width= 0.4pt,dash pattern=on 1pt off 3pt ,line join=round,line cap=round] ( 89.35, 61.20) -- ( 89.35,239.88);

\path[draw=drawColor,line width= 0.4pt,dash pattern=on 1pt off 3pt ,line join=round,line cap=round] (128.71, 61.20) -- (128.71,239.88);

\path[draw=drawColor,line width= 0.4pt,dash pattern=on 1pt off 3pt ,line join=round,line cap=round] (168.07, 61.20) -- (168.07,239.88);

\path[draw=drawColor,line width= 0.4pt,dash pattern=on 1pt off 3pt ,line join=round,line cap=round] (207.44, 61.20) -- (207.44,239.88);

\path[draw=drawColor,line width= 0.4pt,dash pattern=on 1pt off 3pt ,line join=round,line cap=round] (246.80, 61.20) -- (246.80,239.88);

\path[draw=drawColor,line width= 0.4pt,dash pattern=on 1pt off 3pt ,line join=round,line cap=round] (286.16, 61.20) -- (286.16,239.88);

\path[draw=drawColor,line width= 0.4pt,dash pattern=on 1pt off 3pt ,line join=round,line cap=round] (325.52, 61.20) -- (325.52,239.88);

\path[draw=drawColor,line width= 0.4pt,dash pattern=on 1pt off 3pt ,line join=round,line cap=round] ( 49.20, 67.82) -- (336.15, 67.82);

\path[draw=drawColor,line width= 0.4pt,dash pattern=on 1pt off 3pt ,line join=round,line cap=round] ( 49.20,109.18) -- (336.15,109.18);

\path[draw=drawColor,line width= 0.4pt,dash pattern=on 1pt off 3pt ,line join=round,line cap=round] ( 49.20,150.54) -- (336.15,150.54);

\path[draw=drawColor,line width= 0.4pt,dash pattern=on 1pt off 3pt ,line join=round,line cap=round] ( 49.20,191.90) -- (336.15,191.90);

\path[draw=drawColor,line width= 0.4pt,dash pattern=on 1pt off 3pt ,line join=round,line cap=round] ( 49.20,233.26) -- (336.15,233.26);
\definecolor{drawColor}{RGB}{0,0,0}

\path[draw=drawColor,line width= 0.8pt,line join=round,line cap=round] ( 59.83,150.54) -- (305.84,150.54);

\path[draw=drawColor,line width= 0.8pt,line join=round,line cap=round] (299.58,146.93) --
	(305.84,150.54) --
	(299.58,154.15);

\path[draw=drawColor,line width= 0.8pt,line join=round,line cap=round] ( 89.35, 67.82) -- ( 89.35,233.26);

\path[draw=drawColor,line width= 0.8pt,line join=round,line cap=round] ( 92.96,227.00) --
	( 89.35,233.26) --
	( 85.74,227.00);

\path[draw=drawColor,line width= 0.8pt,dash pattern=on 1pt off 3pt ,line join=round,line cap=round] ( 89.35,150.54) circle (196.81);

\node[text=drawColor,anchor=base,inner sep=0pt, outer sep=0pt, scale=  0.70] at (315,128.11) {$Re(z_k)$};

\node[text=drawColor,anchor=base,inner sep=0pt, outer sep=0pt, scale=  0.70] at (125,220) {$Im(z_k)$};

\path[draw=drawColor,line width= 0.4pt,line join=round,line cap=round] (286.16,129.86) --
	(286.16,171.22);

\path[draw=drawColor,line width= 0.4pt,line join=round,line cap=round] ( 87.38,212.58) --
	( 91.32,212.58);

\node[text=drawColor,anchor=base,inner sep=0pt, outer sep=0pt, scale=  0.70] at ( 69.67,210.15) {$0.15$};
\definecolor{fillColor}{RGB}{255,0,0}

\path[fill=fillColor] ( 73.49,150.54) circle (  2.25);

\path[fill=fillColor] ( 78.02,112.74) circle (  2.25);

\path[fill=fillColor] ( 78.02,188.34) circle (  2.25);

\path[fill=fillColor] (281.29,150.54) circle (  2.25);

\path[fill=fillColor] (286.16,150.54) circle (  2.25);

\node[text=drawColor,anchor=base,inner sep=0pt, outer sep=0pt, scale=  1.00] at (266.48,168.64) {$z_2$};

\node[text=drawColor,anchor=base,inner sep=0pt, outer sep=0pt, scale=  1.00] at (305.84,168.64) {$z_1$};

\node[text=drawColor,anchor=base,inner sep=0pt, outer sep=0pt, scale=  1.00] at ( 60,160.36) {$z_3$};

\node[text=drawColor,anchor=base,inner sep=0pt, outer sep=0pt, scale=  1.00] at ( 60,185) {$z_4$};

\node[text=drawColor,anchor=base,inner sep=0pt, outer sep=0pt, scale=  1.00] at ( 60,110) {$z_5$};
\end{scope}
\end{tikzpicture}
\end{floatrow}
\vspace{-1cm}
\caption{\label{Figure1-example2}
[Example~\ref{Numerical-example-2}]
Zeroes of  the characteristic polynomial.}
\end{figure}
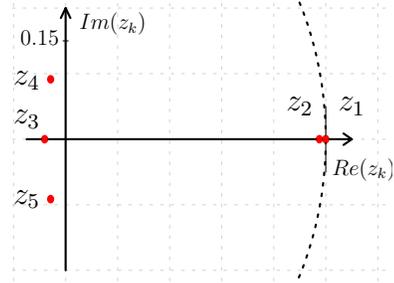

The values of the roots $z_k$ can now be used to set the linear system~\eqref{system2-for-b},
which can be numerically solved to obtain the associated coefficients $b_{k,j}$, see the table in
Figure~\ref{Figure1-example2}. Observe that the coefficients $b_{4,1}$ and $b_{5,1}$
associated with the conjugate roots $z_4$ and $z_5$ are also conjugate.
With these elements we can now calculate
the exact ruin probability $\psi(u)$ and its approximations $\hat{\psi}_1(u)$ and $\hat{\psi}_2(u)$,
which are shown in Figure~\ref{Figure2-example2}. In this case the three probabilities
stay fairly close to each other.
\vspace*{-1.3cm}

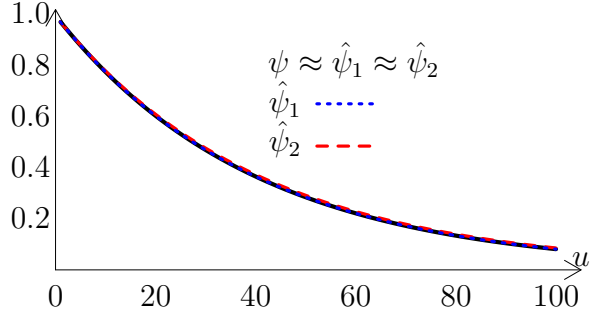
\begin{figure}[H]
\begin{floatrow}
\capbtabbox{%
\footnotesize
%
\begin{tabular}{|c|l|l|l|}
\hline
$u$ & $\psi(u)$ & $\hat{\psi}_1(u)$ & $\hat{\psi}_2(u)$\\ 
\hline
$0$ & $0.99$ & $-$ & $-$ \\
$1$ & $0.9699$ & $0.9704$& $0.9699$  \\
$5$ & $0.8778$ & $0.8778$& $0.8792$  \\
$10$ & $0.7744$ & $0.7744$& $0.7778$  \\
$20$ & $0.6027$ & $0.6027$& $0.6087$  \\
$50$ & $0.2842$ & $0.2842$& $0.2917$  \\
$75$ & $0.1519$ & $0.1519$& $0.1580$  \\
$100$ & $0.0812$ & $0.0812$& $0.0856$  \\
\hline
\end{tabular}
\vspace{2cm}
}{%
}
\begin{tikzpicture}[x=1pt,y=1pt]
\definecolor{fillColor}{RGB}{255,255,255}
\path[use as bounding box,fill=fillColor,fill opacity=0.00] (0,0) rectangle (289.08,216.81);
\begin{scope}
\path[clip] ( 49.20, 61.20) rectangle (263.88,167.61);
\definecolor{drawColor}{RGB}{0,0,0}

\path[draw=drawColor,line width= 1.6pt,line join=round,line cap=round] ( 59.04,158.98) --
	( 60.94,156.70) --
	( 62.83,154.44) --
	( 64.72,152.23) --
	( 66.62,150.07) --
	( 68.51,147.97) --
	( 70.40,145.92) --
	( 72.30,143.92) --
	( 74.19,141.97) --
	( 76.08,140.07) --
	( 77.98,138.21) --
	( 79.87,136.40) --
	( 81.76,134.64) --
	( 83.65,132.92) --
	( 85.55,131.24) --
	( 87.44,129.61) --
	( 89.33,128.01) --
	( 91.23,126.45) --
	( 93.12,124.94) --
	( 95.01,123.46) --
	( 96.91,122.01) --
	( 98.80,120.61) --
	(100.69,119.23) --
	(102.59,117.89) --
	(104.48,116.59) --
	(106.37,115.32) --
	(108.27,114.07) --
	(110.16,112.86) --
	(112.05,111.68) --
	(113.94,110.53) --
	(115.84,109.41) --
	(117.73,108.31) --
	(119.62,107.24) --
	(121.52,106.20) --
	(123.41,105.18) --
	(125.30,104.19) --
	(127.20,103.23) --
	(129.09,102.28) --
	(130.98,101.36) --
	(132.88,100.47) --
	(134.77, 99.59) --
	(136.66, 98.74) --
	(138.56, 97.91) --
	(140.45, 97.10) --
	(142.34, 96.31) --
	(144.23, 95.53) --
	(146.13, 94.78) --
	(148.02, 94.05) --
	(149.91, 93.33) --
	(151.81, 92.64) --
	(153.70, 91.96) --
	(155.59, 91.29) --
	(157.49, 90.64) --
	(159.38, 90.01) --
	(161.27, 89.40) --
	(163.17, 88.80) --
	(165.06, 88.21) --
	(166.95, 87.64) --
	(168.85, 87.08) --
	(170.74, 86.54) --
	(172.63, 86.01) --
	(174.52, 85.49) --
	(176.42, 84.99) --
	(178.31, 84.50) --
	(180.20, 84.02) --
	(182.10, 83.55) --
	(183.99, 83.10) --
	(185.88, 82.65) --
	(187.78, 82.22) --
	(189.67, 81.80) --
	(191.56, 81.38) --
	(193.46, 80.98) --
	(195.35, 80.59) --
	(197.24, 80.21) --
	(199.14, 79.83) --
	(201.03, 79.47) --
	(202.92, 79.12) --
	(204.81, 78.77) --
	(206.71, 78.43) --
	(208.60, 78.10) --
	(210.49, 77.78) --
	(212.39, 77.47) --
	(214.28, 77.17) --
	(216.17, 76.87) --
	(218.07, 76.58) --
	(219.96, 76.29) --
	(221.85, 76.02) --
	(223.75, 75.75) --
	(225.64, 75.49) --
	(227.53, 75.23) --
	(229.43, 74.98) --
	(231.32, 74.74) --
	(233.21, 74.50) --
	(235.10, 74.27) --
	(237.00, 74.04) --
	(238.89, 73.82) --
	(240.78, 73.61) --
	(242.68, 73.40) --
	(244.57, 73.19) --
	(246.46, 72.99);
\definecolor{drawColor}{RGB}{255,0,0}

\path[draw=drawColor,line width= 1.2pt,dash pattern=on 4pt off 4pt ,line join=round,line cap=round] ( 59.04,158.98) --
	( 60.94,156.70) --
	( 62.83,154.49) --
	( 64.72,152.32) --
	( 66.62,150.21) --
	( 68.51,148.15) --
	( 70.40,146.14) --
	( 72.30,144.18) --
	( 74.19,142.26) --
	( 76.08,140.39) --
	( 77.98,138.57) --
	( 79.87,136.79) --
	( 81.76,135.06) --
	( 83.65,133.36) --
	( 85.55,131.71) --
	( 87.44,130.10) --
	( 89.33,128.53) --
	( 91.23,126.99) --
	( 93.12,125.49) --
	( 95.01,124.03) --
	( 96.91,122.60) --
	( 98.80,121.21) --
	(100.69,119.85) --
	(102.59,118.53) --
	(104.48,117.23) --
	(106.37,115.97) --
	(108.27,114.74) --
	(110.16,113.54) --
	(112.05,112.37) --
	(113.94,111.22) --
	(115.84,110.11) --
	(117.73,109.02) --
	(119.62,107.96) --
	(121.52,106.92) --
	(123.41,105.91) --
	(125.30,104.92) --
	(127.20,103.96) --
	(129.09,103.02) --
	(130.98,102.10) --
	(132.88,101.20) --
	(134.77,100.33) --
	(136.66, 99.48) --
	(138.56, 98.65) --
	(140.45, 97.83) --
	(142.34, 97.04) --
	(144.23, 96.27) --
	(146.13, 95.51) --
	(148.02, 94.78) --
	(149.91, 94.06) --
	(151.81, 93.36) --
	(153.70, 92.68) --
	(155.59, 92.01) --
	(157.49, 91.36) --
	(159.38, 90.72) --
	(161.27, 90.10) --
	(163.17, 89.50) --
	(165.06, 88.91) --
	(166.95, 88.33) --
	(168.85, 87.77) --
	(170.74, 87.22) --
	(172.63, 86.69) --
	(174.52, 86.17) --
	(176.42, 85.66) --
	(178.31, 85.16) --
	(180.20, 84.68) --
	(182.10, 84.20) --
	(183.99, 83.74) --
	(185.88, 83.29) --
	(187.78, 82.85) --
	(189.67, 82.42) --
	(191.56, 82.00) --
	(193.46, 81.59) --
	(195.35, 81.20) --
	(197.24, 80.81) --
	(199.14, 80.43) --
	(201.03, 80.06) --
	(202.92, 79.70) --
	(204.81, 79.34) --
	(206.71, 79.00) --
	(208.60, 78.66) --
	(210.49, 78.34) --
	(212.39, 78.02) --
	(214.28, 77.70) --
	(216.17, 77.40) --
	(218.07, 77.10) --
	(219.96, 76.81) --
	(221.85, 76.53) --
	(223.75, 76.26) --
	(225.64, 75.99) --
	(227.53, 75.72) --
	(229.43, 75.47) --
	(231.32, 75.22) --
	(233.21, 74.97) --
	(235.10, 74.73) --
	(237.00, 74.50) --
	(238.89, 74.28) --
	(240.78, 74.05) --
	(242.68, 73.84) --
	(244.57, 73.63) --
	(246.46, 73.42);
\definecolor{drawColor}{RGB}{0,0,255}

\path[draw=drawColor,line width= 1.2pt,dash pattern=on 1pt off 3pt ,line join=round,line cap=round] ( 59.04,159.03) --
	( 60.94,156.70) --
	( 62.83,154.44) --
	( 64.72,152.23) --
	( 66.62,150.07) --
	( 68.51,147.97) --
	( 70.40,145.92) --
	( 72.30,143.92) --
	( 74.19,141.97) --
	( 76.08,140.07) --
	( 77.98,138.21) --
	( 79.87,136.40) --
	( 81.76,134.64) --
	( 83.65,132.92) --
	( 85.55,131.24) --
	( 87.44,129.61) --
	( 89.33,128.01) --
	( 91.23,126.45) --
	( 93.12,124.94) --
	( 95.01,123.46) --
	( 96.91,122.01) --
	( 98.80,120.61) --
	(100.69,119.23) --
	(102.59,117.89) --
	(104.48,116.59) --
	(106.37,115.32) --
	(108.27,114.07) --
	(110.16,112.86) --
	(112.05,111.68) --
	(113.94,110.53) --
	(115.84,109.41) --
	(117.73,108.31) --
	(119.62,107.24) --
	(121.52,106.20) --
	(123.41,105.18) --
	(125.30,104.19) --
	(127.20,103.23) --
	(129.09,102.28) --
	(130.98,101.36) --
	(132.88,100.47) --
	(134.77, 99.59) --
	(136.66, 98.74) --
	(138.56, 97.91) --
	(140.45, 97.10) --
	(142.34, 96.31) --
	(144.23, 95.53) --
	(146.13, 94.78) --
	(148.02, 94.05) --
	(149.91, 93.33) --
	(151.81, 92.64) --
	(153.70, 91.96) --
	(155.59, 91.29) --
	(157.49, 90.64) --
	(159.38, 90.01) --
	(161.27, 89.40) --
	(163.17, 88.80) --
	(165.06, 88.21) --
	(166.95, 87.64) --
	(168.85, 87.08) --
	(170.74, 86.54) --
	(172.63, 86.01) --
	(174.52, 85.49) --
	(176.42, 84.99) --
	(178.31, 84.50) --
	(180.20, 84.02) --
	(182.10, 83.55) --
	(183.99, 83.10) --
	(185.88, 82.65) --
	(187.78, 82.22) --
	(189.67, 81.80) --
	(191.56, 81.38) --
	(193.46, 80.98) --
	(195.35, 80.59) --
	(197.24, 80.21) --
	(199.14, 79.83) --
	(201.03, 79.47) --
	(202.92, 79.12) --
	(204.81, 78.77) --
	(206.71, 78.43) --
	(208.60, 78.10) --
	(210.49, 77.78) --
	(212.39, 77.47) --
	(214.28, 77.17) --
	(216.17, 76.87) --
	(218.07, 76.58) --
	(219.96, 76.29) --
	(221.85, 76.02) --
	(223.75, 75.75) --
	(225.64, 75.49) --
	(227.53, 75.23) --
	(229.43, 74.98) --
	(231.32, 74.74) --
	(233.21, 74.50) --
	(235.10, 74.27) --
	(237.00, 74.04) --
	(238.89, 73.82) --
	(240.78, 73.61) --
	(242.68, 73.40) --
	(244.57, 73.19) --
	(246.46, 72.99);
\definecolor{drawColor}{RGB}{0,0,0}

\path[draw=drawColor,line width= 0.3pt,line join=round,line cap=round] ( 57.15, 65.14) -- (255.93, 65.14);

\path[draw=drawColor,line width= 0.3pt,line join=round,line cap=round] (249.67, 61.53) --
	(255.93, 65.14) --
	(249.67, 68.75);

\path[draw=drawColor,line width= 0.3pt,line join=round,line cap=round] ( 57.15, 65.14) -- ( 57.15,163.67);

\path[draw=drawColor,line width= 0.3pt,line join=round,line cap=round] ( 60.76,157.41) --
	( 57.15,163.67) --
	( 53.54,157.41);

\node[text=drawColor,anchor=base,inner sep=0pt, outer sep=0pt, scale=  1.00] at (170,140) {$\psi\approx\hat{\psi}_1\approx\hat{\psi}_2$};

\node[text=drawColor,anchor=base,inner sep=0pt, outer sep=0pt, scale=  1.00] at (144,125) {$\hat{\psi}_1$};
\node[text=drawColor,anchor=base,inner sep=0pt, outer sep=0pt, scale=  1.00] at (144,110) {$\hat{\psi}_2$};

\definecolor{drawColor}{RGB}{0,0,255}
\path[draw=drawColor,line width= 1.2pt,dash pattern=on 1pt off 3pt ,line join=round,line cap=round] (156,128) --
	(178,128);

\definecolor{drawColor}{RGB}{255,0,0}
\path[draw=drawColor,line width= 1.2pt,dash pattern=on 4pt off 4pt ,line join=round,line cap=round] (156,112) --
	(179,112);

\definecolor{drawColor}{RGB}{0,0,0}
\node[text=drawColor,anchor=base,inner sep=0pt, outer sep=0pt, scale=  1.00] at (255.93, 66.51) {$u$};
\end{scope}
\begin{scope}
\path[clip] (  0.00,  0.00) rectangle (289.08,216.81);
\definecolor{drawColor}{RGB}{0,0,0}

\node[text=drawColor,anchor=base,inner sep=0pt, outer sep=0pt, scale=  1.00] at ( 57.15, 51.60) {0};

\node[text=drawColor,anchor=base,inner sep=0pt, outer sep=0pt, scale=  1.00] at ( 95.01, 51.60) {20};

\node[text=drawColor,anchor=base,inner sep=0pt, outer sep=0pt, scale=  1.00] at (132.88, 51.60) {40};

\node[text=drawColor,anchor=base,inner sep=0pt, outer sep=0pt, scale=  1.00] at (170.74, 51.60) {60};

\node[text=drawColor,anchor=base,inner sep=0pt, outer sep=0pt, scale=  1.00] at (208.60, 51.60) {80};

\node[text=drawColor,anchor=base,inner sep=0pt, outer sep=0pt, scale=  1.00] at (246.46, 51.60) {100};

\node[text=drawColor,anchor=base east,inner sep=0pt, outer sep=0pt, scale=  1.00] at ( 55.20, 81.05) {0.2};

\node[text=drawColor,anchor=base east,inner sep=0pt, outer sep=0pt, scale=  1.00] at ( 55.20,100.40) {0.4};

\node[text=drawColor,anchor=base east,inner sep=0pt, outer sep=0pt, scale=  1.00] at ( 55.20,119.75) {0.6};

\node[text=drawColor,anchor=base east,inner sep=0pt, outer sep=0pt, scale=  1.00] at ( 55.20,139.10) {0.8};

\node[text=drawColor,anchor=base east,inner sep=0pt, outer sep=0pt, scale=  1.00] at ( 55.20,158.45) {1.0};
\end{scope}
\end{tikzpicture}
\end{floatrow}
\vspace{-1.5cm}
\caption{\label{Figure2-example2}
[Example~\ref{Numerical-example-2}]
Ruin probabilities and their approximations.}
\end{figure}

} 


\example{
\label{Numerical-example-3}
Suppose that $m=7$ and claims $Y$ follow the distribution given by the two
probabilities $f(0)={7}/{8}$ and $f(7)={1}/{8}$,
with mean $\mu_Y=\psi(0)={7}/{8}$.
The coefficients $\alpha$ are calculated and the characteristic
polynomial~\eqref{Characteristic-polynomial} can be constructed. The seven
roots and their multiplicities are shown in Figure~\ref{Figure1-example3}.
\vspace*{-0.2cm}

\begin{figure}[H]
\begin{floatrow}
\capbtabbox{%
\footnotesize
%
\begin{tabular}{|c|l|c|l|}
\hline
 & $z_k$ & $n_k$ &  $b_{k,j}$\\ 
\hline
$z_1$ & $1$ & $1$ & $0$ \\
$z_2$ & $0.9577$ & $1$ & $0.9305$\\
$z_3$ & $-0.6556$ & $1$ & $0.0125$\\
$z_4$ & $-0.3674+0.5577i$ & $1$ & $(1.29+0.54i)E-02$\\
$z_5$ & $-0.3674-0.5577i$ & $1$ & $(1.29-0.54i)E-02$\\
$z_6$ & $0.2878+0.6536i$ & $1$ & $(1.56+1.47i)E-02$\\
$z_7$ & $0.2878-0.6536i$ & $1$ & $(1.56-1.47i)E-02$\\
\hline
\end{tabular}
\vspace{1cm}
}{%
}
\begin{tikzpicture}[x=1pt,y=1pt]
\definecolor{fillColor}{RGB}{255,255,255}
\path[use as bounding box,fill=fillColor,fill opacity=0.00] (0,0) rectangle (216.81,144.54);
\begin{scope}
\path[clip] ( 12.00, 12.00) rectangle (204.81,132.54);
\definecolor{fillColor}{RGB}{255,0,0}

\path[fill=fillColor] ( 27.90, 72.27) circle (  1.69);

\path[fill=fillColor] ( 53.58, 35.64) circle (  1.69);

\path[fill=fillColor] ( 53.58,108.90) circle (  1.69);

\path[fill=fillColor] (111.95, 29.34) circle (  1.69);

\path[fill=fillColor] (111.95,115.20) circle (  1.69);

\path[fill=fillColor] (171.63, 72.27) circle (  1.69);

\path[fill=fillColor] (175.40, 72.27) circle (  1.69);
\end{scope}
\begin{scope}
\path[clip] ( 12.00, 12.00) rectangle (204.81,132.54);
\definecolor{drawColor}{RGB}{211,211,211}

\path[draw=drawColor,line width= 0.4pt,dash pattern=on 1pt off 3pt ,line join=round,line cap=round] ( 41.77, 12.00) -- ( 41.77,132.54);

\path[draw=drawColor,line width= 0.4pt,dash pattern=on 1pt off 3pt ,line join=round,line cap=round] ( 86.31, 12.00) -- ( 86.31,132.54);

\path[draw=drawColor,line width= 0.4pt,dash pattern=on 1pt off 3pt ,line join=round,line cap=round] (130.85, 12.00) -- (130.85,132.54);

\path[draw=drawColor,line width= 0.4pt,dash pattern=on 1pt off 3pt ,line join=round,line cap=round] (175.40, 12.00) -- (175.40,132.54);

\path[draw=drawColor,line width= 0.4pt,dash pattern=on 1pt off 3pt ,line join=round,line cap=round] ( 12.00, 39.43) -- (204.81, 39.43);

\path[draw=drawColor,line width= 0.4pt,dash pattern=on 1pt off 3pt ,line join=round,line cap=round] ( 12.00, 72.27) -- (204.81, 72.27);

\path[draw=drawColor,line width= 0.4pt,dash pattern=on 1pt off 3pt ,line join=round,line cap=round] ( 12.00,105.11) -- (204.81,105.11);
\definecolor{drawColor}{RGB}{0,0,0}

\path[draw=drawColor,line width= 0.8pt,line join=round,line cap=round] ( 19.14, 72.27) -- (197.67, 72.27);

\path[draw=drawColor,line width= 0.8pt,line join=round,line cap=round] (191.41, 68.66) --
	(197.67, 72.27) --
	(191.41, 75.88);

\path[draw=drawColor,line width= 0.8pt,line join=round,line cap=round] ( 86.31, 16.46) -- ( 86.31,128.08);

\path[draw=drawColor,line width= 0.8pt,line join=round,line cap=round] ( 89.92,121.82) --
	( 86.31,128.08) --
	( 82.70,121.82);

\path[draw=drawColor,line width= 0.8pt,dash pattern=on 1pt off 3pt ,line join=round,line cap=round] ( 86.31, 72.27) circle ( 89.09);

\node[text=drawColor,anchor=base,inner sep=0pt, outer sep=0pt, scale=  0.70] at (188.76, 60.67) {$Re(z_k)$};

\node[text=drawColor,anchor=base,inner sep=0pt, outer sep=0pt, scale=  0.70] at (104.13,126.33) {$Im(z_k)$};

\path[draw=drawColor,line width= 0.8pt,line join=round,line cap=round] (175.40, 68.99) --
	(175.40, 75.55);

\path[draw=drawColor,line width= 0.4pt,line join=round,line cap=round] ( 85.42,116.91) --
	( 87.20,116.91);

\path[draw=drawColor,line width= 0.4pt,line join=round,line cap=round] ( 85.42, 27.63) --
	( 87.20, 27.63);

\node[text=drawColor,anchor=base,inner sep=0pt, outer sep=0pt, scale=  0.70] at ( 72.95,114.48) {$0.68$};

\node[text=drawColor,anchor=base,inner sep=0pt, outer sep=0pt, scale=  0.70] at ( 72.95, 25.20) {$-0.68$};
\definecolor{fillColor}{RGB}{255,0,0}

\path[fill=fillColor] ( 27.90, 72.27) circle (  1.69);

\path[fill=fillColor] ( 53.58, 35.64) circle (  1.69);

\path[fill=fillColor] ( 53.58,108.90) circle (  1.69);

\path[fill=fillColor] (111.95, 29.34) circle (  1.69);

\path[fill=fillColor] (111.95,115.20) circle (  1.69);

\path[fill=fillColor] (171.63, 72.27) circle (  1.69);

\path[fill=fillColor] (175.40, 72.27) circle (  1.69);

\node[text=drawColor,anchor=base,inner sep=0pt, outer sep=0pt, scale=  1.00] at (184.31, 82.82) {$z_1$};

\node[text=drawColor,anchor=base,inner sep=0pt, outer sep=0pt, scale=  1.00] at (166.49, 82.82) {$z_2$};

\node[text=drawColor,anchor=base,inner sep=0pt, outer sep=0pt, scale=  1.00] at ( 27.87, 82.82) {$z_3$};

\node[text=drawColor,anchor=base,inner sep=0pt, outer sep=0pt, scale=  1.00] at ( 53.62,115.66) {$z_4$};

\node[text=drawColor,anchor=base,inner sep=0pt, outer sep=0pt, scale=  1.00] at ( 53.62, 23.71) {$z_5$};

\node[text=drawColor,anchor=base,inner sep=0pt, outer sep=0pt, scale=  1.00] at (113.04,104.49) {$z_6$};

\node[text=drawColor,anchor=base,inner sep=0pt, outer sep=0pt, scale=  1.00] at (113.04, 34.88) {$z_7$};
\end{scope}
\end{tikzpicture}
\end{floatrow}
\vspace{0cm}
\caption{\label{Figure1-example3}
[Example~\ref{Numerical-example-3}]
Zeroes of  the characteristic polynomial.}
\end{figure}
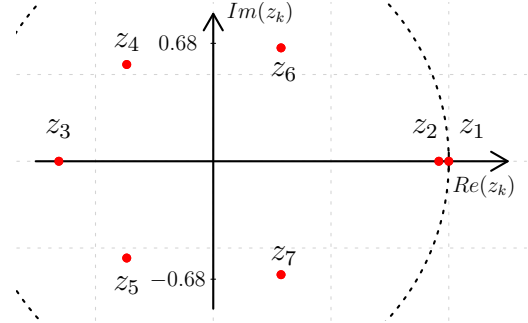

The values of the roots $z_k$ can now be used to set the linear system~\eqref{system2-for-b},
which can be numerically solved to obtain the associated coefficients $b_{k,j}$, see the table in
Figure~\ref{Figure1-example3}.
With these elements we can now calculate the exact
ruin probability $\psi(u)$ and its approximations $\hat{\psi}_1(u)$ and $\hat{\psi}_2(u)$,
which are shown in Figure~\ref{Figure2-example3}.
In this case, the approximation $\hat{\psi}_2(u)$ overestimate the exact ruin probability.

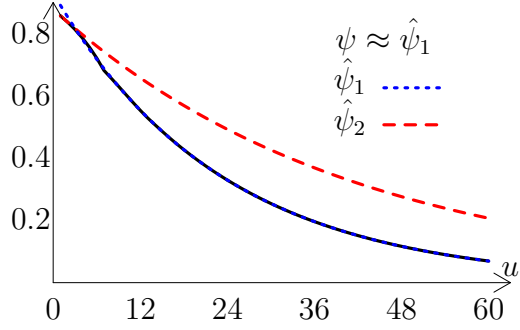
\begin{figure}[H]
\begin{floatrow}
\capbtabbox{%
\footnotesize
%
\begin{tabular}{|c|l|l|l|}
\hline
$u$ & $\psi(u)$ & $\hat{\psi}_1(u)$ & $\hat{\psi}_2(u)$ \\ 
\hline
$0$ & $0.875$ & $-$ & $-$\\
$1$ & $0.8571$ & $0.8911$& $0.8571$ \\
$12$ & $0.5535$ & $0.5537$& $0.6576$ \\
$24$ & $0.3294$ & $0.3294$& $0.4924$ \\
$36$ & $0.1960$ & $0.1960$& $0.3688$ \\
$48$ & $0.1166$ & $0.1166$& $0.2762$ \\
$60$ & $0.0694$ & $0.0694$& $0.2068$ \\
\hline
\end{tabular}
\vspace{0.5cm}
}{%
}
\begin{tikzpicture}[x=1pt,y=1pt]
\definecolor{fillColor}{RGB}{255,255,255}
\path[use as bounding box,fill=fillColor,fill opacity=0.00] (0,0) rectangle (216.81,144.54);
\begin{scope}
\path[clip] ( 18.00, 18.00) rectangle (204.81,132.54);
\definecolor{drawColor}{RGB}{0,0,0}

\path[draw=drawColor,line width= 1.2pt,line join=round,line cap=round] ( 27.66,123.25) --
	( 30.41,120.84) --
	( 33.16,118.09) --
	( 35.90,114.95) --
	( 38.65,111.36) --
	( 41.39,107.26) --
	( 44.14,102.57) --
	( 46.88, 99.62) --
	( 49.63, 96.59) --
	( 52.37, 93.51) --
	( 55.12, 90.45) --
	( 57.87, 87.46) --
	( 60.61, 84.64) --
	( 63.36, 82.07) --
	( 66.10, 79.57) --
	( 68.85, 77.13) --
	( 71.59, 74.80) --
	( 74.34, 72.56) --
	( 77.09, 70.43) --
	( 79.83, 68.40) --
	( 82.58, 66.45) --
	( 85.32, 64.57) --
	( 88.07, 62.78) --
	( 90.81, 61.06) --
	( 93.56, 59.42) --
	( 96.30, 57.85) --
	( 99.05, 56.34) --
	(101.80, 54.89) --
	(104.54, 53.51) --
	(107.29, 52.19) --
	(110.03, 50.92) --
	(112.78, 49.71) --
	(115.52, 48.54) --
	(118.27, 47.43) --
	(121.01, 46.36) --
	(123.76, 45.34) --
	(126.51, 44.36) --
	(129.25, 43.43) --
	(132.00, 42.53) --
	(134.74, 41.67) --
	(137.49, 40.85) --
	(140.23, 40.06) --
	(142.98, 39.31) --
	(145.72, 38.58) --
	(148.47, 37.89) --
	(151.22, 37.23) --
	(153.96, 36.59) --
	(156.71, 35.99) --
	(159.45, 35.40) --
	(162.20, 34.85) --
	(164.94, 34.31) --
	(167.69, 33.80) --
	(170.44, 33.31) --
	(173.18, 32.84) --
	(175.93, 32.40) --
	(178.67, 31.97) --
	(181.42, 31.55) --
	(184.16, 31.16) --
	(186.91, 30.78) --
	(189.65, 30.42);
\definecolor{drawColor}{RGB}{255,0,0}

\path[draw=drawColor,line width= 1.2pt,dash pattern=on 4pt off 4pt ,line join=round,line cap=round] ( 27.66,123.25) --
	( 30.41,120.84) --
	( 33.16,118.49) --
	( 35.90,116.20) --
	( 38.65,113.97) --
	( 41.39,111.78) --
	( 44.14,109.65) --
	( 46.88,107.57) --
	( 49.63,105.54) --
	( 52.37,103.55) --
	( 55.12,101.62) --
	( 57.87, 99.73) --
	( 60.61, 97.88) --
	( 63.36, 96.08) --
	( 66.10, 94.32) --
	( 68.85, 92.61) --
	( 71.59, 90.93) --
	( 74.34, 89.30) --
	( 77.09, 87.70) --
	( 79.83, 86.14) --
	( 82.58, 84.62) --
	( 85.32, 83.14) --
	( 88.07, 81.69) --
	( 90.81, 80.27) --
	( 93.56, 78.89) --
	( 96.30, 77.54) --
	( 99.05, 76.22) --
	(101.80, 74.94) --
	(104.54, 73.68) --
	(107.29, 72.46) --
	(110.03, 71.26) --
	(112.78, 70.10) --
	(115.52, 68.96) --
	(118.27, 67.84) --
	(121.01, 66.76) --
	(123.76, 65.70) --
	(126.51, 64.66) --
	(129.25, 63.65) --
	(132.00, 62.67) --
	(134.74, 61.71) --
	(137.49, 60.77) --
	(140.23, 59.85) --
	(142.98, 58.95) --
	(145.72, 58.08) --
	(148.47, 57.23) --
	(151.22, 56.39) --
	(153.96, 55.58) --
	(156.71, 54.79) --
	(159.45, 54.01) --
	(162.20, 53.25) --
	(164.94, 52.52) --
	(167.69, 51.80) --
	(170.44, 51.09) --
	(173.18, 50.41) --
	(175.93, 49.73) --
	(178.67, 49.08) --
	(181.42, 48.44) --
	(184.16, 47.82) --
	(186.91, 47.21) --
	(189.65, 46.61);
\definecolor{drawColor}{RGB}{0,0,255}

\path[draw=drawColor,line width= 1.2pt,dash pattern=on 1pt off 3pt ,line join=round,line cap=round] ( 27.66,127.25) --
	( 30.41,122.80) --
	( 33.16,118.54) --
	( 35.90,114.47) --
	( 38.65,110.56) --
	( 41.39,106.82) --
	( 44.14,103.24) --
	( 46.88, 99.81) --
	( 49.63, 96.53) --
	( 52.37, 93.38) --
	( 55.12, 90.37) --
	( 57.87, 87.48) --
	( 60.61, 84.72) --
	( 63.36, 82.08) --
	( 66.10, 79.54) --
	( 68.85, 77.12) --
	( 71.59, 74.79) --
	( 74.34, 72.57) --
	( 77.09, 70.44) --
	( 79.83, 68.40) --
	( 82.58, 66.44) --
	( 85.32, 64.57) --
	( 88.07, 62.78) --
	( 90.81, 61.06) --
	( 93.56, 59.42) --
	( 96.30, 57.84) --
	( 99.05, 56.34) --
	(101.80, 54.89) --
	(104.54, 53.51) --
	(107.29, 52.19) --
	(110.03, 50.92) --
	(112.78, 49.71) --
	(115.52, 48.54) --
	(118.27, 47.43) --
	(121.01, 46.36) --
	(123.76, 45.34) --
	(126.51, 44.36) --
	(129.25, 43.43) --
	(132.00, 42.53) --
	(134.74, 41.67) --
	(137.49, 40.85) --
	(140.23, 40.06) --
	(142.98, 39.31) --
	(145.72, 38.58) --
	(148.47, 37.89) --
	(151.22, 37.23) --
	(153.96, 36.59) --
	(156.71, 35.99) --
	(159.45, 35.40) --
	(162.20, 34.85) --
	(164.94, 34.31) --
	(167.69, 33.80) --
	(170.44, 33.31) --
	(173.18, 32.84) --
	(175.93, 32.40) --
	(178.67, 31.97) --
	(181.42, 31.55) --
	(184.16, 31.16) --
	(186.91, 30.78) --
	(189.65, 30.42);
\definecolor{drawColor}{RGB}{0,0,0}

\path[draw=drawColor,line width= 0.3pt,line join=round,line cap=round] ( 24.92, 22.24) -- (197.89, 22.24);

\path[draw=drawColor,line width= 0.3pt,line join=round,line cap=round] (191.63, 18.63) --
	(197.89, 22.24) --
	(191.63, 25.86);

\path[draw=drawColor,line width= 0.3pt,line join=round,line cap=round] ( 24.92, 22.24) -- ( 24.92,128.30);

\path[draw=drawColor,line width= 0.3pt,line join=round,line cap=round] ( 28.53,122.04) --
	( 24.92,128.30) --
	( 21.31,122.04);

\node[text=drawColor,anchor=base,inner sep=0pt, outer sep=0pt, scale=  1.00] at (150,110.75) {$\psi\approx\hat{\psi}_1$};

\node[text=drawColor,anchor=base,inner sep=0pt, outer sep=0pt, scale=  1.00] at (137.29,95.60) {$\hat{\psi}_1$};
\node[text=drawColor,anchor=base,inner sep=0pt, outer sep=0pt, scale=  1.00] at (137.29,80) {$\hat{\psi}_2$};

\definecolor{drawColor}{RGB}{0,0,255}
\path[draw=drawColor,line width= 1.2pt,dash pattern=on 1pt off 3pt ,line join=round,line cap=round] (150.01,96) --
	(172,96);

\definecolor{drawColor}{RGB}{255,0,0}
\path[draw=drawColor,line width= 1.2pt,dash pattern=on 4pt off 4pt ,line join=round,line cap=round] (150,82) --
	(173,82);
	
\definecolor{drawColor}{RGB}{0,0,0}
\node[text=drawColor,anchor=base,inner sep=0pt, outer sep=0pt, scale=  1.00] at (197.89, 24.66) {$u$};
\end{scope}
\begin{scope}
\path[clip] (  0.00,  0.00) rectangle (216.81,144.54);
\definecolor{drawColor}{RGB}{0,0,0}

\node[text=drawColor,anchor=base,inner sep=0pt, outer sep=0pt, scale=  1.00] at ( 24.92,  8.40) {0};

\node[text=drawColor,anchor=base,inner sep=0pt, outer sep=0pt, scale=  1.00] at ( 57.87,  8.40) {12};

\node[text=drawColor,anchor=base,inner sep=0pt, outer sep=0pt, scale=  1.00] at ( 90.81,  8.40) {24};

\node[text=drawColor,anchor=base,inner sep=0pt, outer sep=0pt, scale=  1.00] at (123.76,  8.40) {36};

\node[text=drawColor,anchor=base,inner sep=0pt, outer sep=0pt, scale=  1.00] at (156.71,  8.40) {48};

\node[text=drawColor,anchor=base,inner sep=0pt, outer sep=0pt, scale=  1.00] at (189.65,  8.40) {60};

\node[text=drawColor,anchor=base east,inner sep=0pt, outer sep=0pt, scale=  1.00] at ( 24.00, 42.37) {0.2};

\node[text=drawColor,anchor=base east,inner sep=0pt, outer sep=0pt, scale=  1.00] at ( 24.00, 65.93) {0.4};

\node[text=drawColor,anchor=base east,inner sep=0pt, outer sep=0pt, scale=  1.00] at ( 24.00, 89.50) {0.6};

\node[text=drawColor,anchor=base east,inner sep=0pt, outer sep=0pt, scale=  1.00] at ( 24.00,113.07) {0.8};
\end{scope}
\end{tikzpicture}
\end{floatrow}
\vspace{-0.2cm}
\caption{\label{Figure2-example3}
[Example~\ref{Numerical-example-3}]
Ruin probabilities and their approximations.}
\end{figure}

} 


\example{
\label{Numerical-example-4}
Suppose that $m=7$ and claims $Y$ follow the distribution given by
\begin{displaymath}
f(0)=\frac{1}{2},\, f(1)=\frac{3}{7},\, f(2)=\frac{3}{392},\, f(3)=\frac{145}{2744},\, f(4)=\frac{775}{76832}, \end{displaymath}
\begin{displaymath} f(5)=\frac{219}{268912},\, f(6)=\frac{67}{2151296},\, f(7)=\frac{1}{2151296},
\end{displaymath}
with mean $\mu_Y=\psi(0)=0.647$.
The coefficients $\alpha$ are calculated and the characteristic
polynomial~\eqref{Characteristic-polynomial} can be constructed. The
roots and their multiplicities are shown in Figure~\ref{Figure1-example4}. The roots in this example are all real,
we have the two positive roots $z_1$ and $z_2$, and a negative root $z_3$ with multiplicity $5$.

\begin{figure}[H]
\begin{floatrow}
\capbtabbox{%
\footnotesize
%
\begin{tabular}{|c|c|c|c|}
\hline
 & $z_k$ & $n_k$ &  $b_{k,j}$\\ 
\hline
$z_1$ & $1$ & $1$ & $0$ \\
$z_2$ & $1/2$ & $1$ & $0.7242$\\
$z_3$ & $-1/14$ & $5$ & $0.2758$\\
      &         &     & $0.4150$\\
      &         &     & $0.2133$\\
      &         &     & $0.0454$\\
      &         &     & $0.0034$\\                  
\hline
\end{tabular}
\vspace{1cm}
}{%
}
\begin{tikzpicture}[x=1pt,y=1pt]
\definecolor{fillColor}{RGB}{255,255,255}
\path[use as bounding box,fill=fillColor,fill opacity=0.00] (0,0) rectangle (216.81,144.54);
\begin{scope}
\path[clip] ( 12.00, 12.00) rectangle (204.81,132.54);
\definecolor{fillColor}{RGB}{255,0,0}

\path[fill=fillColor] ( 41.13, 72.27) circle (  1.69);

\path[fill=fillColor] (111.48, 72.27) circle (  1.69);

\path[fill=fillColor] (173.04, 72.27) circle (  1.69);
\end{scope}
\begin{scope}
\path[clip] ( 12.00, 12.00) rectangle (204.81,132.54);
\definecolor{drawColor}{RGB}{211,211,211}

\path[draw=drawColor,line width= 0.4pt,dash pattern=on 1pt off 3pt ,line join=round,line cap=round] ( 49.92, 12.00) -- ( 49.92,132.54);

\path[draw=drawColor,line width= 0.4pt,dash pattern=on 1pt off 3pt ,line join=round,line cap=round] (111.48, 12.00) -- (111.48,132.54);

\path[draw=drawColor,line width= 0.4pt,dash pattern=on 1pt off 3pt ,line join=round,line cap=round] (173.04, 12.00) -- (173.04,132.54);

\path[draw=drawColor,line width= 0.4pt,dash pattern=on 1pt off 3pt ,line join=round,line cap=round] ( 12.00, 16.46) -- (204.81, 16.46);

\path[draw=drawColor,line width= 0.4pt,dash pattern=on 1pt off 3pt ,line join=round,line cap=round] ( 12.00, 44.37) -- (204.81, 44.37);

\path[draw=drawColor,line width= 0.4pt,dash pattern=on 1pt off 3pt ,line join=round,line cap=round] ( 12.00, 72.27) -- (204.81, 72.27);

\path[draw=drawColor,line width= 0.4pt,dash pattern=on 1pt off 3pt ,line join=round,line cap=round] ( 12.00,100.17) -- (204.81,100.17);

\path[draw=drawColor,line width= 0.4pt,dash pattern=on 1pt off 3pt ,line join=round,line cap=round] ( 12.00,128.08) -- (204.81,128.08);
\definecolor{drawColor}{RGB}{0,0,0}

\path[draw=drawColor,line width= 0.4pt,line join=round,line cap=round] ( 19.14, 72.27) -- (185.36, 72.27);

\path[draw=drawColor,line width= 0.4pt,line join=round,line cap=round] (179.10, 68.66) --
	(185.36, 72.27) --
	(179.10, 75.88);

\path[draw=drawColor,line width= 0.4pt,line join=round,line cap=round] ( 49.92, 16.46) -- ( 49.92,128.08);

\path[draw=drawColor,line width= 0.4pt,line join=round,line cap=round] ( 53.54,121.82) --
	( 49.92,128.08) --
	( 46.31,121.82);

\path[draw=drawColor,line width= 0.8pt,dash pattern=on 1pt off 3pt ,line join=round,line cap=round] ( 49.92, 72.27) circle (123.12);

\node[text=drawColor,anchor=base,inner sep=0pt, outer sep=0pt, scale=  0.75] at (185.36, 56.44) {$Re(z_k)$};

\node[text=drawColor,anchor=base,inner sep=0pt, outer sep=0pt, scale=  0.75] at ( 68.39,115.04) {$Im(z_k)$};

\path[draw=drawColor,line width= 0.4pt,line join=round,line cap=round] (173.04, 58.32) --
	(173.04, 86.22);

\path[draw=drawColor,line width= 0.4pt,line join=round,line cap=round] ( 48.69,114.12) --
	( 51.15,114.12);

\node[text=drawColor,anchor=base,inner sep=0pt, outer sep=0pt, scale=  0.75] at ( 37.61,111.52) {$0.15$};
\definecolor{fillColor}{RGB}{255,0,0}

\path[fill=fillColor] ( 41.13, 72.27) circle (  1.69);

\path[fill=fillColor] (111.48, 72.27) circle (  1.69);

\path[fill=fillColor] (173.04, 72.27) circle (  1.69);

\node[text=drawColor,anchor=base,inner sep=0pt, outer sep=0pt, scale=  1.00] at (111.48, 83.64) {$z_2$};

\node[text=drawColor,anchor=base,inner sep=0pt, outer sep=0pt, scale=  1.00] at ( 41.13, 83.64) {$z_3$};

\node[text=drawColor,anchor=base,inner sep=0pt, outer sep=0pt, scale=  1.00] at (185.36, 83.64) {$z_1$};
\end{scope}
\end{tikzpicture}
\end{floatrow}
\vspace{-0.3cm}
\caption{\label{Figure1-example4}
[Example~\ref{Numerical-example-4}]
Zeroes of  the characteristic polynomial.}
\end{figure}
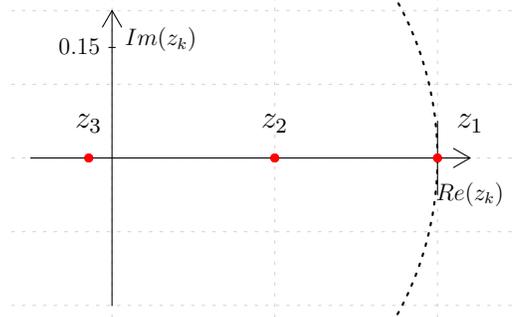

The values of the roots $z_k$ can now be used to set the linear system~\eqref{system2-for-b},
which can be numerically solved to obtain the associated coefficients $b_{k,j}$.
With these elements we can now calculate
the exact ruin probability $\psi(u)$ and its approximations $\hat{\psi}_1(u)$ and $\hat{\psi}_2(u)$,
which are shown in Figure~\ref{Figure2-example4}.
Observe that $\hat{\psi}_2(u)=\psi(u)$ for
$u=1,2$, but slightly overestimates $\psi(u)$ for $u\ge3$.

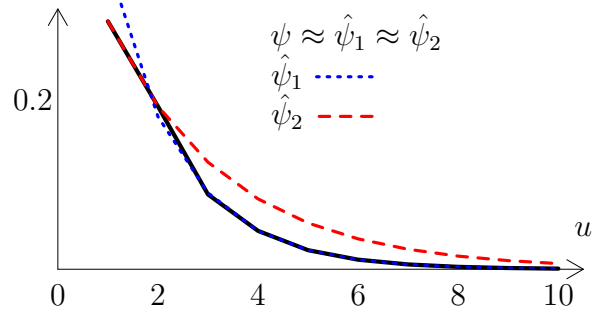
\begin{figure}[H]
\begin{floatrow}
\capbtabbox{%
\footnotesize
%
\begin{tabular}{|c|c|c|c|}
\hline
$u$ & $\psi(u)$ & $\hat{\psi}_1(u)$ & $\hat{\psi}_2(u)$ \\ 
\hline
$0$  & $0.6470$ & $-$ & $-$\\
$1$  & $0.2940$ & $0.3621$& $0.2940$ \\
$2$  & $0.1932$ & $0.1810$& $0.1932$ \\
$4$  & $0.0455$ & $0.0453$& $0.0834$ \\
$6$  & $0.0113$ & $0.0113$& $0.0360$ \\
$8$  & $0.0028$ & $0.0028$& $0.0155$ \\
$10$ & $0.0007$ & $0.0007$& $0.0067$ \\
\hline
\end{tabular}
\vspace{2.5cm}
}{%
}
\begin{tikzpicture}[x=1pt,y=1pt]
\definecolor{fillColor}{RGB}{255,255,255}
\path[use as bounding box,fill=fillColor,fill opacity=0.00] (0,0) rectangle (289.08,216.81);
\begin{scope}
\path[clip] ( 49.20, 61.20) rectangle (263.88,167.61);
\definecolor{drawColor}{RGB}{0,0,0}

\path[draw=drawColor,line width= 1.6pt,line join=round,line cap=round] ( 76.08,158.98) --
	( 95.01,126.79) --
	(113.94, 93.46) --
	(132.88, 79.66) --
	(151.81, 72.35) --
	(170.74, 68.75) --
	(189.67, 66.95) --
	(208.60, 66.04) --
	(227.53, 65.59) --
	(246.46, 65.37);
\definecolor{drawColor}{RGB}{255,0,0}

\path[draw=drawColor,line width= 1.2pt,dash pattern=on 4pt off 4pt ,line join=round,line cap=round] ( 76.08,158.98) --
	( 95.01,126.79) --
	(113.94,105.65) --
	(132.88, 91.75) --
	(151.81, 82.62) --
	(170.74, 76.63) --
	(189.67, 72.69) --
	(208.60, 70.10) --
	(227.53, 68.40) --
	(246.46, 67.28);
\definecolor{drawColor}{RGB}{0,0,255}

\path[draw=drawColor,line width= 1.2pt,dash pattern=on 1pt off 3pt ,line join=round,line cap=round] ( 76.08,180.70) --
	( 95.01,122.92) --
	(113.94, 94.03) --
	(132.88, 79.59) --
	(151.81, 72.36) --
	(170.74, 68.75) --
	(189.67, 66.95) --
	(208.60, 66.04) --
	(227.53, 65.59) --
	(246.46, 65.37);
\definecolor{drawColor}{RGB}{0,0,0}

\path[draw=drawColor,line width= 0.3pt,line join=round,line cap=round] ( 57.15, 65.14) -- (255.93, 65.14);

\path[draw=drawColor,line width= 0.3pt,line join=round,line cap=round] (249.67, 61.53) --
	(255.93, 65.14) --
	(249.67, 68.75);

\path[draw=drawColor,line width= 0.3pt,line join=round,line cap=round] ( 57.15, 65.14) -- ( 57.15,163.67);

\path[draw=drawColor,line width= 0.3pt,line join=round,line cap=round] ( 60.76,157.41) --
	( 57.15,163.67) --
	( 53.54,157.41);

\node[text=drawColor,anchor=base,inner sep=0pt, outer sep=0pt, scale=  1.00] at (170,150) {$\psi\approx\hat{\psi}_1\approx\hat{\psi}_2$};

\node[text=drawColor,anchor=base,inner sep=0pt, outer sep=0pt, scale=  1.00] at (144,135) {$\hat{\psi}_1$};
\node[text=drawColor,anchor=base,inner sep=0pt, outer sep=0pt, scale=  1.00] at (144,120) {$\hat{\psi}_2$};

\definecolor{drawColor}{RGB}{0,0,255}
\path[draw=drawColor,line width= 1.2pt,dash pattern=on 1pt off 3pt ,line join=round,line cap=round] (155.01,138) --
	(177,138);

\definecolor{drawColor}{RGB}{255,0,0}
\path[draw=drawColor,line width= 1.2pt,dash pattern=on 4pt off 4pt ,line join=round,line cap=round] (156,123) --
	(179,123);

\definecolor{drawColor}{RGB}{0,0,0}
\node[text=drawColor,anchor=base,inner sep=0pt, outer sep=0pt, scale=  1.00] at (255.93, 77.63) {$u$};
\end{scope}
\begin{scope}
\path[clip] (  0.00,  0.00) rectangle (289.08,216.81);
\definecolor{drawColor}{RGB}{0,0,0}

\node[text=drawColor,anchor=base,inner sep=0pt, outer sep=0pt, scale=  1.00] at ( 57.15, 51.60) {0};

\node[text=drawColor,anchor=base,inner sep=0pt, outer sep=0pt, scale=  1.00] at ( 95.01, 51.60) {2};

\node[text=drawColor,anchor=base,inner sep=0pt, outer sep=0pt, scale=  1.00] at (132.88, 51.60) {4};

\node[text=drawColor,anchor=base,inner sep=0pt, outer sep=0pt, scale=  1.00] at (170.74, 51.60) {6};

\node[text=drawColor,anchor=base,inner sep=0pt, outer sep=0pt, scale=  1.00] at (208.60, 51.60) {8};

\node[text=drawColor,anchor=base,inner sep=0pt, outer sep=0pt, scale=  1.00] at (246.46, 51.60) {10};

\node[text=drawColor,anchor=base east,inner sep=0pt, outer sep=0pt, scale=  1.00] at ( 55.20,125.52) {0.2};
\end{scope}
\end{tikzpicture}
\end{floatrow}
\vspace{-1.5cm}
\caption{\label{Figure2-example4}
[Example~\ref{Numerical-example-4}]
Ruin probabilities and their approximations.}
\end{figure}

} 


\example{
\label{Numerical-example-5}
Suppose that $m=8$ and claims $Y$ follow the distribution given by
\begin{displaymath}
f(0)=\frac{1}{2},\, f(1)=\frac{9}{28},\, f(2)=\frac{477}{3136},\, f(3)=\frac{543}{21952},\, f(4)=\frac{9433}{19668992}, \end{displaymath}
\begin{displaymath} f(5)=\frac{4462}{3813049},\, f(6)=\frac{146689}{1927561216},\, f(7)=\frac{7155}{1927561216}, f(8)=\frac{2809}{1927561216},
\end{displaymath}
with mean $\mu_Y=\psi(0)=0.70810867$.
The coefficients $\alpha$ are calculated and the characteristic
polynomial~\eqref{Characteristic-polynomial} can be constructed. The
roots and their multiplicities are shown in Figure~\ref{Figure1-example5}.

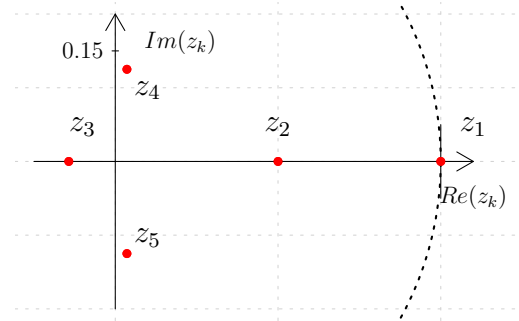
\begin{figure}[H]
\begin{floatrow}
\capbtabbox{%
\footnotesize
%
\begin{tabular}{|c|l|c|l|}
\hline
 & $z_k$ & $n_k$ &  $b_{k,j}$\\ 
\hline
$z_1$ & $1$ & $1$ & $0$ \\
$z_2$ & $1/2$ & $1$ & $0.82594$\\
$z_3$ & $-1/7$ & $2$ & $0.07341$\\
      &         &     & $0.02094$\\
$z_4$ & $1/28+i/8$ & $2$ & $0.01243-0.00945i$\\
      &            &     & $0.05033-0.03952i$\\
$z_5$ & $1/28-i/8$ & $2$ & $0.01243+0.00945i$\\
      &            &     & $0.05033+0.03952i$\\      
\hline
\end{tabular}
\vspace{0.5cm}
}{%
}
\begin{tikzpicture}[x=1pt,y=1pt]
\definecolor{fillColor}{RGB}{255,255,255}
\path[use as bounding box,fill=fillColor,fill opacity=0.00] (0,0) rectangle (216.81,144.54);
\begin{scope}
\path[clip] ( 12.00, 12.00) rectangle (204.81,132.54);
\definecolor{fillColor}{RGB}{255,0,0}

\path[fill=fillColor] ( 54.32, 37.39) circle (  1.69);

\path[fill=fillColor] ( 54.32,107.15) circle (  1.69);

\path[fill=fillColor] ( 32.33, 72.27) circle (  1.69);

\path[fill=fillColor] (111.48, 72.27) circle (  1.69);

\path[fill=fillColor] (173.04, 72.27) circle (  1.69);
\end{scope}
\begin{scope}
\path[clip] ( 12.00, 12.00) rectangle (204.81,132.54);
\definecolor{drawColor}{RGB}{211,211,211}

\path[draw=drawColor,line width= 0.4pt,dash pattern=on 1pt off 3pt ,line join=round,line cap=round] ( 49.92, 12.00) -- ( 49.92,132.54);

\path[draw=drawColor,line width= 0.4pt,dash pattern=on 1pt off 3pt ,line join=round,line cap=round] (111.48, 12.00) -- (111.48,132.54);

\path[draw=drawColor,line width= 0.4pt,dash pattern=on 1pt off 3pt ,line join=round,line cap=round] (173.04, 12.00) -- (173.04,132.54);

\path[draw=drawColor,line width= 0.4pt,dash pattern=on 1pt off 3pt ,line join=round,line cap=round] ( 12.00, 16.46) -- (204.81, 16.46);

\path[draw=drawColor,line width= 0.4pt,dash pattern=on 1pt off 3pt ,line join=round,line cap=round] ( 12.00, 44.37) -- (204.81, 44.37);

\path[draw=drawColor,line width= 0.4pt,dash pattern=on 1pt off 3pt ,line join=round,line cap=round] ( 12.00, 72.27) -- (204.81, 72.27);

\path[draw=drawColor,line width= 0.4pt,dash pattern=on 1pt off 3pt ,line join=round,line cap=round] ( 12.00,100.17) -- (204.81,100.17);

\path[draw=drawColor,line width= 0.4pt,dash pattern=on 1pt off 3pt ,line join=round,line cap=round] ( 12.00,128.08) -- (204.81,128.08);
\definecolor{drawColor}{RGB}{0,0,0}

\path[draw=drawColor,line width= 0.4pt,line join=round,line cap=round] ( 19.14, 72.27) -- (185.36, 72.27);

\path[draw=drawColor,line width= 0.4pt,line join=round,line cap=round] (179.10, 68.66) --
	(185.36, 72.27) --
	(179.10, 75.88);

\path[draw=drawColor,line width= 0.4pt,line join=round,line cap=round] ( 49.92, 16.46) -- ( 49.92,128.08);

\path[draw=drawColor,line width= 0.4pt,line join=round,line cap=round] ( 53.54,121.82) --
	( 49.92,128.08) --
	( 46.31,121.82);

\path[draw=drawColor,line width= 0.8pt,dash pattern=on 1pt off 3pt ,line join=round,line cap=round] ( 49.92, 72.27) circle (123.12);

\node[text=drawColor,anchor=base,inner sep=0pt, outer sep=0pt, scale=  0.75] at (185.36, 56.44) {$Re(z_k)$};

\node[text=drawColor,anchor=base,inner sep=0pt, outer sep=0pt, scale=  0.75] at ( 74.55,115.04) {$Im(z_k)$};

\path[draw=drawColor,line width= 0.4pt,line join=round,line cap=round] (173.04, 58.32) --
	(173.04, 86.22);

\path[draw=drawColor,line width= 0.4pt,line join=round,line cap=round] ( 48.69,114.12) --
	( 51.15,114.12);

\node[text=drawColor,anchor=base,inner sep=0pt, outer sep=0pt, scale=  0.75] at ( 37.61,111.52) {$0.15$};
\definecolor{fillColor}{RGB}{255,0,0}

\path[fill=fillColor] ( 54.32, 37.39) circle (  1.69);

\path[fill=fillColor] ( 54.32,107.15) circle (  1.69);

\path[fill=fillColor] ( 32.33, 72.27) circle (  1.69);

\path[fill=fillColor] (111.48, 72.27) circle (  1.69);

\path[fill=fillColor] (173.04, 72.27) circle (  1.69);

\node[text=drawColor,anchor=base,inner sep=0pt, outer sep=0pt, scale=  1.00] at (111.48, 83.64) {$z_2$};

\node[text=drawColor,anchor=base,inner sep=0pt, outer sep=0pt, scale=  1.00] at ( 37.61, 83.64) {$z_3$};

\node[text=drawColor,anchor=base,inner sep=0pt, outer sep=0pt, scale=  1.00] at ( 62.23, 97.59) {$z_4$};

\node[text=drawColor,anchor=base,inner sep=0pt, outer sep=0pt, scale=  1.00] at ( 62.23, 41.78) {$z_5$};

\node[text=drawColor,anchor=base,inner sep=0pt, outer sep=0pt, scale=  1.00] at (185.36, 83.64) {$z_1$};
\end{scope}
\end{tikzpicture}
\end{floatrow}
\vspace{-0.5cm}
\caption{\label{Figure1-example5}
[Example~\ref{Numerical-example-5}]
Zeroes of  the characteristic polynomial.}
\end{figure}

The values of the roots $z_k$ can now be used to set the linear system~\eqref{system2-for-b},
which can be numerically solved to obtain the associated coefficients $b_{k,j}$.
With these elements we can now calculate
the ruin probability $\psi(u)$ and its approximations $\hat{\psi}_1(u)$ and $\hat{\psi}_2(u)$,
which are shown in Figure~\ref{Figure2-example5}. The two approximations seem to be
very accurate.
\vspace*{-1.5cm}

\begin{figure}[H]
\begin{floatrow}
\capbtabbox{%
\footnotesize
%
\begin{tabular}{|c|l|l|l|}
\hline
$u$ & $\psi(u)$ & $\hat{\psi}_1(u)$ & $\hat{\psi}_2(u)$ \\ 
\hline
$0$ & $0.7081$ & $-$ & $-$\\
$1$ & $0.4162$ & $0.4130$& $0.4162$ \\
$2$ & $0.2077$ & $0.2065$& $0.2077$ \\
$4$ & $0.0517$ & $0.0516$& $0.0517$ \\
$6$ & $0.0129$ & $0.0129$ & $0.0129$ \\
$8$ & $0.0032$ & $0.0032$ & $0.0032$\\
$10$ & $0.0008$ & $0.0008$ & $0.0008$\\
\hline
\end{tabular}
\vspace{2.2cm}
}{%
}
\begin{tikzpicture}[x=1pt,y=1pt]
\definecolor{fillColor}{RGB}{255,255,255}
\path[use as bounding box,fill=fillColor,fill opacity=0.00] (0,0) rectangle (289.08,216.81);
\begin{scope}
\path[clip] ( 24.00, 48.00) rectangle (277.08,168.81);
\definecolor{drawColor}{RGB}{0,0,0}

\path[draw=drawColor,line width= 1.6pt,line join=round,line cap=round] ( 55.69,154.17) --
	( 78.01,103.23) --
	(100.33, 77.48) --
	(122.64, 65.10) --
	(144.96, 58.78) --
	(167.28, 55.63) --
	(189.60, 54.05) --
	(211.91, 53.26) --
	(234.23, 52.87) --
	(256.55, 52.67);
\definecolor{drawColor}{RGB}{255,0,0}

\path[draw=drawColor,line width= 1.2pt,dash pattern=on 4pt off 4pt ,line join=round,line cap=round] ( 55.69,154.17) --
	( 78.01,103.23) --
	(100.33, 77.80) --
	(122.64, 65.12) --
	(144.96, 58.78) --
	(167.28, 55.62) --
	(189.60, 54.05) --
	(211.91, 53.26) --
	(234.23, 52.87) --
	(256.55, 52.67);
\definecolor{drawColor}{RGB}{0,0,255}

\path[draw=drawColor,line width= 1.2pt,dash pattern=on 1pt off 3pt ,line join=round,line cap=round] ( 55.69,153.37) --
	( 78.01,102.92) --
	(100.33, 77.70) --
	(122.64, 65.09) --
	(144.96, 58.78) --
	(167.28, 55.63) --
	(189.60, 54.05) --
	(211.91, 53.26) --
	(234.23, 52.87) --
	(256.55, 52.67);
\definecolor{drawColor}{RGB}{0,0,0}

\path[draw=drawColor,line width= 0.3pt,line join=round,line cap=round] ( 33.37, 52.47) -- (267.71, 52.47);

\path[draw=drawColor,line width= 0.3pt,line join=round,line cap=round] (261.45, 48.86) --
	(267.71, 52.47) --
	(261.45, 56.09);

\path[draw=drawColor,line width= 0.3pt,line join=round,line cap=round] ( 33.37, 52.47) -- ( 33.37,164.34);

\path[draw=drawColor,line width= 0.3pt,line join=round,line cap=round] ( 36.99,158.08) --
	( 33.37,164.34) --
	( 29.76,158.08);

\node[text=drawColor,anchor=base,inner sep=0pt, outer sep=0pt, scale=  1.00] at (164,140) {$\psi\approx\hat{\psi}_1\approx\hat{\psi}_2$};

\node[text=drawColor,anchor=base,inner sep=0pt, outer sep=0pt, scale=  1.00] at (137.29,125) {$\hat{\psi}_1$};
\node[text=drawColor,anchor=base,inner sep=0pt, outer sep=0pt, scale=  1.00] at (137.29,110) {$\hat{\psi}_2$};

\definecolor{drawColor}{RGB}{0,0,255}
\path[draw=drawColor,line width= 1.2pt,dash pattern=on 1pt off 3pt ,line join=round,line cap=round] (150.01,128) --
	(172,128);

\definecolor{drawColor}{RGB}{255,0,0}
\path[draw=drawColor,line width= 1.2pt,dash pattern=on 4pt off 4pt ,line join=round,line cap=round] (151,113) --
	(174,113);

\definecolor{drawColor}{RGB}{0,0,0}
\node[text=drawColor,anchor=base,inner sep=0pt, outer sep=0pt, scale=  1.00] at (267.71, 61.22) {$u$};
\end{scope}
\begin{scope}
\path[clip] (  0.00,  0.00) rectangle (289.08,216.81);
\definecolor{drawColor}{RGB}{0,0,0}

\node[text=drawColor,anchor=base,inner sep=0pt, outer sep=0pt, scale=  1.00] at ( 33.37, 38.40) {0};

\node[text=drawColor,anchor=base,inner sep=0pt, outer sep=0pt, scale=  1.00] at ( 78.01, 38.40) {2};

\node[text=drawColor,anchor=base,inner sep=0pt, outer sep=0pt, scale=  1.00] at (122.64, 38.40) {4};

\node[text=drawColor,anchor=base,inner sep=0pt, outer sep=0pt, scale=  1.00] at (167.28, 38.40) {6};

\node[text=drawColor,anchor=base,inner sep=0pt, outer sep=0pt, scale=  1.00] at (211.91, 38.40) {8};

\node[text=drawColor,anchor=base,inner sep=0pt, outer sep=0pt, scale=  1.00] at (256.55, 38.40) {10};

\node[text=drawColor,anchor=base east,inner sep=0pt, outer sep=0pt, scale=  1.00] at ( 27.60, 73.46) {0.1};

\node[text=drawColor,anchor=base east,inner sep=0pt, outer sep=0pt, scale=  1.00] at ( 27.60, 97.90) {0.2};

\node[text=drawColor,anchor=base east,inner sep=0pt, outer sep=0pt, scale=  1.00] at ( 27.60,122.33) {0.3};

\node[text=drawColor,anchor=base east,inner sep=0pt, outer sep=0pt, scale=  1.00] at ( 27.60,146.76) {0.4};
\end{scope}
\end{tikzpicture}
\end{floatrow}
\vspace{-1cm}
\caption{\label{Figure2-example5}
[Example~\ref{Numerical-example-5}]
Ruin probabilities and their approximations.}
\end{figure}
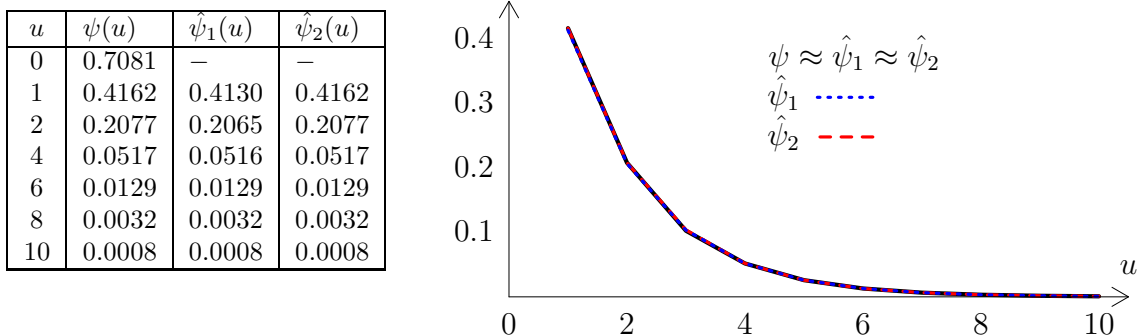

} 


\section{Conclusions}
\label{conclusions}

The theory of linear recurrence sequences has been applied to obtain a new
formula~\eqref{Solution-psi(u)} for the ultimate ruin probability in a discrete-time risk
process with claims following an arbitrary discrete distribution with bounded support.
The formula found is expressed
in terms of the roots of the characteristic polynomial~\eqref{Characteristic-polynomial}
associated with the recurrence sequence~\eqref{Recursive-relation-4}.
Thus, a ruin probability problem has been translated into the classical mathematical problem
of finding the zeroes of a polynomial and that of solving a system of linear equations.
In the particular case when all the roots are simple the reduced
formula~\eqref{Solution-for-simple-roots} is obtained.  We have also presented a simple but
effective
approximation formula~\eqref{Approximation-1} based on the numerical observation that
$b_2z_2^u$ is the leading term in the formula~\eqref{Solution-psi(u)}.\\

Several numerical examples were shown where the roots were found using the {\tt R} pracma package. For these particular examples the exact value of the ruin
probability was calculated along with the approximations proposed. However, no attempts
were made to measure the accuracy of the approximations.\\


As is well known, many results in the theory of ruin involve recursive relations and it is very
likely that the theory of recurrence sequences can be applied in more instances than the one
shown in this work.\\

{\bf Acknowledgments.} We are grateful for the comments and suggestions from anonymous
reviewers and editors. Their corrections helped improve the presentation of our paper.\\

{\bf Funding details.} The authors received no financial support for the research, authorship,
and/or publication of this article.\\

{\bf Data availability statement.}
All data generated or analyzed during this study are included in this published article.\\


{\bf Declaration of competing interest.} The authors declare that they have no known competing
financial interests or personal relationships that could have appeared to influence the work
reported in this paper.\\

\printbibliography

\end{document}